\documentclass[12pt, eqno]{article}
\usepackage{amssymb,latexsym}
\usepackage{amsmath,latexsym}
\usepackage{amscd,latexsym} 
\usepackage[all]{xy}

\begin{document}

\title
{Minimal submanifolds in certain types of kaehler product manifold}

\author{Xingda Liu
\hspace{10mm}Bang Xiao
\thanks{Corresponding author.Supported by the NNSF 11526040 of China.
Partially supported by the NNSF 11371386 of China}\\
{\normalsize Department of Mathematics,  Chongqing University of Technology},\\
  {\normalsize Chongqing 400054, The People's Republic
 of China}\\
{\normalsize (xingdaliu0924@cqut.edu.cn\hspace{5mm}bangxiao@cqut.edu.cn)}}
\maketitle \vspace{-0.1in}

\abstract{
Let $M$ be a real $l$-dimensional minimal submanifold with flat normal connection
in a kaehler product manifold $\overline{M}^m\times \overline{M}^n$ where
$\overline{M}^m$ and $\overline{M}^n$ are complex $m$-dimensional and
complex $n$-dimensional kaehler manifolds with constant holomorphic sectional
curvature $c_1$ and $c_2$ respectively.
We give a formula for the Laplacian of the second fundamental form of $M$. Specially
we discuss the F-anti invariant case.
We also give some applications of this formula.
}\\

{\bf Key words:}\,Kaehler product manifold,
minimal submanifolds.\\

{\bf 2000 Mathematics Subject Classification:}\, 58A05; 53C55.

\section{Introduction}

The complex projective space $CP^m$ and its submanifolds have been studied by many researchers.
For the product manifold K.Yano and
M.Kon \cite{YK1} give the general results.
The F-invariant submanifolds, invariant submanifolds, totally real submanifolds,
of Kaehler product manifold have been studied, for example \cite{MA}, \cite{MSK}, \cite{XN},
\cite{AC}.

In this paper we give a formula for the Laplacian of the second fundamentai form of
a $l$-dimensional minimal submanifold with flat normal connection in
a kaehler product manifold $\overline{M}^m\times \overline{M}^n$, where
$\overline{M}^m$ and $\overline{M}^n$ are complex $m$-dimensional and
complex $n$-dimensional kaehler manifolds with constant holomorphic sectional
curvature $c_1$ and $c_2$ respectively.
We also give some applications of this formula.
Specially in the F-anti invariant case,
we get several simpler results.

\section{Preliminaries}
Let $\overline{M}^m$ be a Kaehlerian manifold of complex dimension $m$
with constant holomorphic sectional curvature $c_1$ and $\overline{M}^n$
be a Kaehlerian manifold of complex dimension $n$ with constant holomorphic sectional
curvature $c_2$. We denote by $J_m$ and $J_n$ almost complex structures of
$\overline{M}^m$ and $\overline{M}^n$ respectively.
We consider the kaeherian product $\overline{M}=\overline{M}^m\times\overline{M}^n$
with product metric $g$ and put
$$JX=J_m\overline{P}X+J_n\overline{Q}X$$
for any vector field $X$ on $\overline{M}$, where $\overline{P}$ and $\overline{Q}$
denote the projection operators. Then we have
$$J_m\overline{P}=\overline{P}J,~~~J_n\overline{Q}=\overline{Q}J,~~~FJ=JF,~~~F^2=I,$$
$$J^2=-I,~~~g(JX,JY)=g(X,Y),~~~\overline{\nabla}_XJ=0,$$
where $F=\overline{P}-\overline{Q}$ is a almost product structure on $\overline{M}$.
Then the Riemannian curvature tensor $\overline{R}$ is given by
\begin{eqnarray}
\overline{R}(X,Y)Z&=&\frac{c_1+c_2}{16}[g(Y,Z)X-g(X,Z)Y+g(JY,Z)JX \nonumber\\
&&-g(JX,Z)JY+2g(X,JY)JZ+2g(FY,Z)FX \nonumber\\
&&-g(FX,Z)FY+g(FJY,Z)FJX \nonumber\\
&&-g(FJX,Z)FJY+2g(FX,JY)FJZ] \\
&&+\frac{c_1-c_2}{16}[g(FY,Z)X-g(FX,Z)Y+g(Y,Z)FX \nonumber\\
&&-g(X,Z)FY+g(FJY,Z)JX-g(FJX,Z)JY \nonumber\\
&&+g(JY,Z)FJX-g(JX,Z)FJY \nonumber\\
&&+2g(FX,JY)JZ+2g(X,JY)JFZ] \nonumber
\end{eqnarray}
for any vector fields $X,~Y$ and $Z$ on $\overline{M}$.

Let $M$ be a real $l$-dimensional submanifold of $\overline{M}$ with induced metric.
We denote by $\overline{\nabla}$ the Levi-Civita connection in $\overline{M}$ and by
$\nabla$ the connection induced on $M$.
Then the Gauss and Weingarten formulas are given respectively by
$$\overline{\nabla}_XY=\nabla_XY+B(X,Y),~~\overline{\nabla}_XV=-A_VX+D_XV$$
for any vector field $X$ and $Y$ tangent to $M$ and any vector field $V$ normal to $M$,
where $D$ denote the normal connection.

For any vector field $X$ tangent to $M$ we put
$$JX=\mathcal{P}X+\mathcal{F}X,~~~~~FX=fX+hX,$$
where $\mathcal{P}X$ and $\mathcal{F}X$ are the the tangential part and normal part of $JX$ respectively,
$fX$ and $hX$ are the tangential part and normal part of $FX$ respectively.

For any vector field $V$ normal to $M$ we put
$$JV=\mathcal{T}V+\mathcal{N}V,~~~~~FV=tV+sV,$$
where $\mathcal{T}V$ and $\mathcal{N}V$ are the tangential part
and normal part of $JV$ respectively, $tV$ and $sV$ are the
tangential part and normal part of $FV$ respectively.
Then $F$, $f$, $s$,  $ht$ and $\mathcal{F}\mathcal{T}$ are symmetric with respect to $g$,~
$\mathcal{P}$,~$\mathcal{N}$,~$\mathcal{F}J$,~$f\mathcal{P}+t\mathcal{F}$ and $\mathcal{F}t+\mathcal{N}s$ are skew-symmetric with respect to $g$,
and $g(\mathcal{F}X,V)=-g(X,\mathcal{T}V)$,~$g(hX,V)=g(X,tV)$,~$g((h\mathcal{P}+s\mathcal{F})X,V)=-g(X,(\mathcal{P}t+\mathcal{T}s)V)$.
We also have
\begin{eqnarray}
\mathcal{P}^2=-I-\mathcal{T}\mathcal{F},~~\mathcal{F}\mathcal{P}+\mathcal{N}\mathcal{F}=0,
~~\mathcal{P}\mathcal{T}+\mathcal{T}\mathcal{N}=0,~~\mathcal{N}^2=-I-\mathcal{F}\mathcal{T}.
\end{eqnarray}
\begin{eqnarray}f^2=I-th,~~hf+sh=0,~~ft+ts=0£¬~~s^2=I-ht,
\end{eqnarray}
\begin{eqnarray}
f\mathcal{T}+t\mathcal{N}=\mathcal{P}t+\mathcal{T}s,~~~h\mathcal{P}+s\mathcal{F}=\mathcal{F}f+\mathcal{N}h,
\end{eqnarray}
\begin{eqnarray}
f\mathcal{P}+t\mathcal{F}=\mathcal{P}f+\mathcal{T}h,~~~h\mathcal{T}+s\mathcal{N}=\mathcal{F}t+\mathcal{N}s,
\end{eqnarray}
\begin{eqnarray}
FJX=JFX&=&(f\mathcal{P}+t\mathcal{F})X+(h\mathcal{P}+s\mathcal{F})X \nonumber\\
&=&(\mathcal{P}f+\mathcal{T}h)X+(\mathcal{F}f+\mathcal{N}h)X,
\end{eqnarray}
\begin{eqnarray}
FJV=JFV&=&(\mathcal{P}t+\mathcal{T}s)V+(\mathcal{F}t+\mathcal{N}s)V\nonumber\\
&=&(f\mathcal{T}+t\mathcal{N})V+(h\mathcal{T}+s\mathcal{N})V.
\end{eqnarray}

Additionally, we have the following relations:
$$tr(f\mathcal{P})=tr(s\mathcal{N})=tr(h\mathcal{T})=tr(\mathcal{F}t)=0,~~tr[(h\mathcal{T})^2]=tr[(\mathcal{F}t)^2]$$

Now we define the covariant derivatives of $\mathcal{P}$,~$\mathcal{F}$,~$\mathcal{T}$ and $\mathcal{N}$
respectively by
$$(\nabla_X\mathcal{P})Y=\nabla_X(\mathcal{P}Y)-\mathcal{P}(\nabla_XY),~
(\nabla_X\mathcal{F})Y=D_X(\mathcal{F}Y)-\mathcal{F}(\nabla_XY),$$
$$(\nabla_X\mathcal{T})V=\nabla_X(\mathcal{T}V)-\mathcal{T}(D_XV),~
(\nabla_X\mathcal{N})V=D_X(\mathcal{N}V)-\mathcal{N}(D_XV).$$
Then we have
$$(\nabla_X\mathcal{P})Y=A_{\mathcal{F}Y}X+\mathcal{T}B(X,Y),~
(\nabla_X\mathcal{F})Y=-B(X,\mathcal{P}Y)+\mathcal{N}B(X,Y),$$
$$(\nabla_X\mathcal{T})V=-\mathcal{P}A_VX+A_{\mathcal{N}V}X,~
(\nabla_X\mathcal{N})V=-\mathcal{F}A_VX-B(X,\mathcal{T}V).$$

Because $\overline{M}$ is equipped with product metric we know that
$\overline{\nabla}F=\overline{\nabla}~\overline{P}=\overline{\nabla}~\overline{Q}=0$ (cf \cite{XN}).
Next we define the covariant derivatives of $f$,~$h$,~$t$, and $s$ by
$$(\nabla_Xf)Y=\nabla_X(fY)-f(\nabla_XY),~~(\nabla_Xh)Y=D_X(hY)-h(\nabla_XY),$$
$$(\nabla_Xt)V=\nabla_X(tV)-t(D_XV),~~(\nabla_Xs)V=D_X(sV)-s(D_XV).$$
Then we have
$$(\nabla_Xf)Y=A_{hY}X+tB(x,y),~~(\nabla_Xh)Y=sB(X,Y)-B(X,fY),$$
$$(\nabla_Xt)V=A_{sV}X-fA_VX,~~(\nabla_Xs)V=-B(X,tV)-hA_VX$$

Denote by $R$ the Riemannian curvature tensor of $M$,
then the equation of Gauss is given by
$$R(X,Y)Z=(\overline{R}(X,Y)Z)^\top+A_{B(Y,Z)}X-A_{B(X,Z)}Y,$$
that is
\begin{eqnarray}
R(X,Y)Z&=&\frac{c_1+c_2}{16}[g(Y,Z)X-g(X,Z)Y+g(\mathcal{P}Y,Z)\mathcal{P}X-g(\mathcal{P}X,Z)\mathcal{P}Y \nonumber\\
&&+2g(X,\mathcal{P}Y)\mathcal{P}Z+2g(fY,Z)fX-g(fX,Z)fY \nonumber\\
&&+g((f\mathcal{P}+t\mathcal{F})Y,Z)(f\mathcal{P}+t\mathcal{F})X-g((f\mathcal{P}+t\mathcal{F})X,Z)(f\mathcal{P}+t\mathcal{F})Y \nonumber\\
&&+2g(X,(f\mathcal{P}+t\mathcal{F})Y)(f\mathcal{P}+t\mathcal{F})Z]   \nonumber\\
&&+\frac{c_1-c_2}{16}[g(fY,Z)X-g(fX,Z)Y+g(Y,Z)fX-g(X,Z)fY \nonumber\\
&&+g((f\mathcal{P}+t\mathcal{F})Y,Z)\mathcal{P}X-g((f\mathcal{P}+t\mathcal{F})X,Z)\mathcal{P}Y  \\
&&+g(\mathcal{P}Y,Z)(f\mathcal{P}+t\mathcal{F})X-g(\mathcal{P}X,Z)(f\mathcal{P}+t\mathcal{F})Y \nonumber\\
&&+2g(X,(f\mathcal{P}+t\mathcal{F})Y)\mathcal{P}Z+2g(X,\mathcal{P}Y)(f\mathcal{P}+t\mathcal{F})Z] \nonumber\\
&&+A_{B(Y,Z)}X-A_{B(X,Z)}Y. \nonumber
\end{eqnarray}

The equation of Codazzi is given by
$$(\nabla_XB)(Y,Z)-(\nabla_YB)(X,Z)=(\overline{R}(X,Y)Z)^\bot,$$
that is
\begin{eqnarray}
&&(\nabla_XB)(Y,Z)-(\nabla_YB)(X,Z) \nonumber\\
&=&\frac{c_1+c_2}{16}[g(\mathcal{P}Y,Z)\mathcal{F}X-g(\mathcal{P}X,Z)\mathcal{F}Y+2g(X,\mathcal{P}Y)\mathcal{F}Z+2g(fY,Z)hX \nonumber\\
&&-g(fX,Z)hY+g((f\mathcal{P}+t\mathcal{F})Y,Z)(h\mathcal{P}+s\mathcal{F})X \nonumber\\
&&-g((f\mathcal{P}+t\mathcal{F})X,Z)(h\mathcal{P}+s\mathcal{F})Y+2g(X,(f\mathcal{P}+t\mathcal{F})Y)(h\mathcal{P}+s\mathcal{F})Z] \nonumber\\
&&+\frac{c_1-c_2}{16}[g(Y,Z)hX-g(X,Z)hY+g((f\mathcal{P}+t\mathcal{F})Y,Z)\mathcal{F}X  \nonumber\\
&&-g((f\mathcal{P}+t\mathcal{F})X,Z)\mathcal{F}Y+g(\mathcal{P}Y,Z)(h\mathcal{P}+s\mathcal{F})X \nonumber\\
&&-g(\mathcal{P}X,Z)(h\mathcal{P}+s\mathcal{F})Y+2g(X,(f\mathcal{P}+t\mathcal{F})Y)\mathcal{F}Z \nonumber\\
&&+2g(X,\mathcal{P}Y)(h\mathcal{P}+s\mathcal{F})Z ]
\end{eqnarray}

The equation of Ricci is given by
$$R^\bot(X,Y)U=(\overline{R}(X,Y)U)^\bot+B(X,A_UY)-B(Y,A_UX)$$
or equivalently
$$g(R^\bot(X,Y)U,V)+g([A_V,A_U]X,Y)=g((\overline{R}(X,Y)U)^\bot,V),$$
where $[A_V,A_U]=A_VA_U-A_UA_V$.
That is
\begin{eqnarray}
&&g(R^\bot(X,Y)U,V)+g([A_V,A_U]X,Y) \nonumber\\
&=&\frac{c_1+c_2}{16}[g(\mathcal{F}Y,U)g(\mathcal{F}X,V)-g(\mathcal{F}X,U)g(\mathcal{F}Y,V) \nonumber\\
&&+2g(X,\mathcal{P}Y)g(\mathcal{N}U,V)+2g(hY,U)g(hX,V)-g(hX,U)g(hY,V) \nonumber\\
&&+g((h\mathcal{P}+s\mathcal{F})Y,U)g((h\mathcal{P}+s\mathcal{F})X,V) \nonumber\\
&&-g((h\mathcal{P}+s\mathcal{F})X,U)g((h\mathcal{P}+s\mathcal{F})Y,V)  \\
&&+2g(X,(f\mathcal{P}+t\mathcal{F})Y)((\mathcal{F}t+\mathcal{N}s)U,V) ] \nonumber\\
&&+\frac{c_1-c_2}{16}[g((h\mathcal{P}+s\mathcal{F})Y,U)g(\mathcal{F}X,V)-g((h\mathcal{P}+s\mathcal{F})X,U)g(\mathcal{F}Y,V)\nonumber\\
&&+g(\mathcal{F}Y,U)g((h\mathcal{P}+s\mathcal{F})X,V)-g(\mathcal{F}X,U)g((h\mathcal{P}+s\mathcal{F})Y,V) \nonumber\\
&&+2g(X,(f\mathcal{P}+t\mathcal{F})Y)g(\mathcal{N}U,V)+2g(X,\mathcal{P}Y)g((\mathcal{F}t+\mathcal{N}s)U,V)] \nonumber
\end{eqnarray}

In the following, we denote by $A_\alpha$ the second fundamental form in the direction of $v_\alpha$,
where $\{v_1,\cdots, v_p\}$ is an orthonormal basis for $T_x(M)^\bot$, $P=2m+2n-l$. We denote by $|\cdot|$
the length of the tensor.
~~~~~~~~~~~

~~~~~~~~~~~~

\textbf{Lemma 2.1}~~Let $M$ be an $l$-dimensional submanifold in $\overline{M}^m\times\overline{M}^n$.
If the normal connection of $M$ is flat, then
\begin{eqnarray}
&&\sum_{\alpha,\beta}g([A_\beta, A_\alpha]\mathcal{T}v_\beta, \mathcal{T}v_\alpha) \nonumber\\
&=&\sum_{\alpha,\beta}[g(A_\alpha\mathcal{T}v_\beta, A_\beta\mathcal{T}v_\alpha)
-g(A_\alpha\mathcal{T}v_\alpha, A_\beta\mathcal{T}v_\beta)] \nonumber\\
&=& \frac{c_1+c_2}{16}\Big\{[tr(\mathcal{F}\mathcal{T})]^2-tr[(\mathcal{F}\mathcal{T})^2]
+2tr(\mathcal{N}\mathcal{F}\mathcal{P}\mathcal{T})-tr[(\mathcal{F}t)^2]  \\
&&+[tr(h\mathcal{P}+s\mathcal{F})\mathcal{T}]^2-tr[(h\mathcal{P}+s\mathcal{F})\mathcal{T}]^2
+2tr[(\mathcal{F}t+\mathcal{N}s)\mathcal{F}(f\mathcal{P}+t\mathcal{F})\mathcal{T}] \Big\}  \nonumber\\
&&+2\frac{c_1-c_2}{16}\Big\{tr(\mathcal{F}\mathcal{T})tr((h\mathcal{P}+s\mathcal{F})\mathcal{T})
-tr(\mathcal{F}\mathcal{T}(h\mathcal{P}+s\mathcal{F})\mathcal{T})  \nonumber\\
&&+tr(\mathcal{N}\mathcal{F}(f\mathcal{P}+t\mathcal{F})\mathcal{T})+tr(\mathcal{F}\mathcal{P}\mathcal{T}(\mathcal{F}t+\mathcal{N}s))  \Big\} \nonumber
\end{eqnarray}
where
\begin{eqnarray}
-tr[(\mathcal{F}\mathcal{T})^2]=-\sum_\alpha g(\mathcal{F}\mathcal{T}v_\alpha,\mathcal{F}\mathcal{T}v_\alpha)\leq 0
\end{eqnarray}
\begin{eqnarray}
&&\sum_{i,\alpha}g([A_{\mathcal{N}v_\alpha},A_\alpha]e_i,\mathcal{P}e_i)=2\sum_{\alpha}tr(A_\alpha A_{\mathcal{N}v_\alpha}\mathcal{P}) \nonumber\\
&=& \frac{c_1+c_2}{16}\Big\{-2tr(\mathcal{T}\mathcal{N}\mathcal{F}\mathcal{P})
-2tr(\mathcal{P}^2)tr(\mathcal{N}^2)+3tr(\mathcal{N}h\mathcal{P}t)  \\
&&-2tr[(\mathcal{P}t+\mathcal{T}s)\mathcal{N}(h\mathcal{P}+s\mathcal{F})\mathcal{P}]
-2tr[(f\mathcal{P}+t\mathcal{F})\mathcal{P}]\cdot tr[(\mathcal{F}t+\mathcal{N}s)\mathcal{N}]\Big\}  \nonumber\\
&&-2\frac{c_1-c_2}{16}\Big\{2tr(\mathcal{P}\mathcal{T}\mathcal{N}(h\mathcal{P}+s\mathcal{F}))
+tr(\mathcal{N}^2)tr((f\mathcal{P}+t\mathcal{F})\mathcal{P})  \nonumber\\
&&+tr(\mathcal{P}^2)tr((\mathcal{F}t+\mathcal{N}s)\mathcal{N}) \Big\}  \nonumber
\end{eqnarray}
\begin{eqnarray}
&&\sum_{\alpha,\beta}[g(A_\alpha tv_\beta, A_\beta tv_\alpha)-g(A_\alpha tv_\alpha, A_\beta tv_\beta)]
=g([A_\alpha,A_\beta]tv_\alpha, tv_\beta)  \\
&=&\frac{c_1+c_2}{16}\Big[(tr(\mathcal{F}t))^2-tr((\mathcal{F}t)^2)-2tr(\mathcal{N}h\mathcal{P}t)+2(tr(ht))^2-tr((ht)^2)\nonumber\\
&&+(tr((h\mathcal{P}+s\mathcal{F})t))^2-tr(((h\mathcal{P}+s\mathcal{F})t)^2)-2tr((\mathcal{F}t+\mathcal{N}s)h(f\mathcal{P}+t\mathcal{F})t)\Big]\nonumber\\
&&+2\frac{c_1-c_2}{16}\Big[-tr((h\mathcal{P}+s\mathcal{F})t\mathcal{F}t)
-tr((f\mathcal{P}+t\mathcal{F})t\mathcal{N}h)-tr(h\mathcal{P}t(\mathcal{F}t+\mathcal{N}s))  \Big]           \nonumber
\end{eqnarray}
where
\begin{eqnarray}
-tr((ht)^2)=-\sum_\alpha g(htv_\alpha,htv_\alpha)\leq 0
\end{eqnarray}
\begin{eqnarray}
&&g([A_{(\mathcal{F}t+\mathcal{N}s)v_\alpha},A_\alpha]e_i,(f\mathcal{P}+t\mathcal{F})e_i) \nonumber\\
&=&2tr(A_\alpha A_{(\mathcal{F}t+\mathcal{N}s)v_\alpha}(f\mathcal{P}+t\mathcal{F})) \\
&=&\frac{c_1+c_2}{16}\Big\{-2tr((f\mathcal{P}+t\mathcal{F})\mathcal{T}(\mathcal{F}t+\mathcal{N}s)\mathcal{F})
-2tr(\mathcal{P}(f\mathcal{P}+t\mathcal{F}))tr(\mathcal{N}(\mathcal{F}t+\mathcal{N}s))  \nonumber\\
&&+3tr(t(\mathcal{F}t+\mathcal{N}s)h(f\mathcal{P}+t\mathcal{F}))
-2tr((\mathcal{P}t+\mathcal{T}s)(\mathcal{F}t+\mathcal{N}s)(h\mathcal{P}+s\mathcal{F})(f\mathcal{P}+t\mathcal{F})) \nonumber\\
&&-2tr((f\mathcal{P}+t\mathcal{F})^2)\cdot tr((\mathcal{F}t+\mathcal{N}s)^2) \Big\} \nonumber\\
&&-2\frac{c_1-c_2}{16}\Big\{tr(\mathcal{T}(\mathcal{F}t+\mathcal{N}s)(h\mathcal{P}+s\mathcal{F})(f\mathcal{P}+t\mathcal{F})) \nonumber\\
&&+tr((\mathcal{P}t+\mathcal{T}s)(\mathcal{F}t+\mathcal{N}s)\mathcal{F}(f\mathcal{P}+t\mathcal{F})) \nonumber\\
&&+tr((f\mathcal{P}+t\mathcal{F})^2)tr(\mathcal{N}(\mathcal{F}t+\mathcal{N}s))
+tr((\mathcal{F}t+\mathcal{N}s)^2)tr(\mathcal{P}(f\mathcal{P}+t\mathcal{F})) \Big\} \nonumber
\end{eqnarray}
where
\begin{eqnarray}
&&-2tr((f\mathcal{P}+t\mathcal{F})^2)\cdot tr((\mathcal{F}t+\mathcal{N}s)^2)  \nonumber\\
&=&-2\sum_{i,\alpha}g((f\mathcal{P}+t\mathcal{F})e_i,(f\mathcal{P}+t\mathcal{F})e_i)\cdot
g((\mathcal{F}t+\mathcal{N}s)v_\alpha,(\mathcal{F}t+\mathcal{N}s)v_\alpha)  \nonumber\\
&\leq &0
\end{eqnarray}
\begin{eqnarray}
&&\sum_{\alpha,\beta}\Big[g(A_\alpha(\mathcal{P}t+\mathcal{T}s)v_\beta, A_\beta(\mathcal{P}t+\mathcal{T}s)v_\alpha)
-g(A_\alpha(\mathcal{P}t+\mathcal{T}s)v_\alpha, A_\beta(\mathcal{P}t+\mathcal{T}s)v_\beta)\Big] \nonumber\\
&=&\sum_{\alpha,\beta}g([A_\alpha,A_\beta](\mathcal{P}t+\mathcal{T}s)v_\alpha,(\mathcal{P}t+\mathcal{T}s)v_\beta) \\
&=&\frac{c_1+c_2}{16}\Big\{[tr(\mathcal{F}(\mathcal{P}t+\mathcal{T}s))]^2-tr[(\mathcal{F}(\mathcal{P}t+\mathcal{T}s))^2]  \nonumber\\
&&+2tr(\mathcal{N}(h\mathcal{P}+s\mathcal{F})\mathcal{P}(\mathcal{P}t+\mathcal{T}s))+2[tr(h(\mathcal{P}t+\mathcal{T}s))]^2
-tr[(h(\mathcal{P}t+\mathcal{T}s))^2]  \nonumber\\
&&+[tr((h\mathcal{P}+s\mathcal{F})(\mathcal{P}t+\mathcal{T}s))]^2-tr[((h\mathcal{P}+s\mathcal{F})(\mathcal{P}t+\mathcal{T}s))^2]  \nonumber\\
&&+2tr((\mathcal{F}t+\mathcal{N}s)(h\mathcal{P}+s\mathcal{F})(f\mathcal{P}+t\mathcal{F})(\mathcal{P}t+\mathcal{T}s)) \Big\}   \nonumber\\
&&+2\frac{c_1-c_2}{16}\Big\{tr(\mathcal{F}(\mathcal{P}t+\mathcal{T}s))\cdot tr((h\mathcal{P}+s\mathcal{F})(\mathcal{P}t+\mathcal{T}s))  \nonumber\\
&&-tr(\mathcal{F}(\mathcal{P}t+\mathcal{T}s)(h\mathcal{P}+s\mathcal{F})(\mathcal{P}t+\mathcal{T}s))
+tr(\mathcal{N}(h\mathcal{P}+s\mathcal{F})(f\mathcal{P}+t\mathcal{F})(\mathcal{P}t+\mathcal{T}s))  \nonumber\\
&&+tr((\mathcal{F}t+\mathcal{N}s)(h\mathcal{P}+s\mathcal{F})\mathcal{P}(\mathcal{P}t+\mathcal{T}s))  \Big\}\nonumber
\end{eqnarray}
where
\begin{eqnarray}
&&-tr[((h\mathcal{P}+s\mathcal{F})(\mathcal{P}t+\mathcal{T}s))^2]    \nonumber\\
&=&-\sum_\alpha g((h\mathcal{P}+s\mathcal{F})(\mathcal{P}t+\mathcal{T}s)v_\alpha,
(h\mathcal{P}+s\mathcal{F})(\mathcal{P}t+\mathcal{T}s)v_\alpha)   \nonumber\\
&\leq& 0
\end{eqnarray}
\begin{eqnarray}
&&g([A_{(\mathcal{F}t+\mathcal{N}s)v_\alpha},A_\alpha]e_i,\mathcal{P}e_i)=2tr(A_\alpha A_{(\mathcal{F}t+\mathcal{N}s)v_\alpha}\mathcal{P})\\
&=&\frac{c_1+c_2}{16}\Big\{-2tr(\mathcal{F}\mathcal{P}\mathcal{T}(\mathcal{F}t+\mathcal{N}s))
-2tr(\mathcal{P}^2)tr(\mathcal{N}(\mathcal{F}t+\mathcal{N}s))  \nonumber\\
&&+3tr(h\mathcal{P}t(\mathcal{F}t+\mathcal{N}s))
-2tr((\mathcal{P}t+\mathcal{T}s)(\mathcal{F}t+\mathcal{N}s)(h\mathcal{P}+s\mathcal{F})\mathcal{P})  \nonumber\\
&&-2tr((f\mathcal{P}+t\mathcal{F})\mathcal{P})\cdot tr((\mathcal{F}t+\mathcal{N}s)^2)\Big\} \nonumber\\
&&+\frac{c_1-c_2}{16}\Big\{-4tr(\mathcal{P}\mathcal{T}(\mathcal{F}t+\mathcal{N}s)(h\mathcal{P}+s\mathcal{F}))  \nonumber\\
&&-2tr((f\mathcal{P}+t\mathcal{F})\mathcal{P})\cdot tr(\mathcal{N}(\mathcal{F}t+\mathcal{N}s))
-2tr(\mathcal{P}^2)\cdot tr((\mathcal{F}t+\mathcal{N}s)^2)         \Big\}\nonumber
\end{eqnarray}
where
\begin{eqnarray}
&&-2tr(\mathcal{P}^2)\cdot tr((\mathcal{F}t+\mathcal{N}s)^2)  \nonumber\\
&=&-2g(\mathcal{P}e_i,\mathcal{P}e_i)\cdot g((\mathcal{F}t+\mathcal{N}s)v_\alpha,(\mathcal{F}t+\mathcal{N}s)v_\alpha)\leq 0
\end{eqnarray}

\begin{eqnarray}
&&\sum_{\alpha,\beta}\Big[ g(A_\alpha\mathcal{T}v_\beta, A_\beta(\mathcal{P}t+\mathcal{T}s)v_\alpha)
-g(A_\alpha\mathcal{T}v_\alpha, A_\beta(\mathcal{P}t+\mathcal{T}s)v_\beta)\Big]  \nonumber \\
&=&\sum_{\alpha,\beta}g([A_\alpha,A_\beta](\mathcal{P}t+\mathcal{T}s)v_\alpha,\mathcal{T}v_\beta)\\
&=&\frac{c_1+c_2}{16}\Big\{tr(\mathcal{F}\mathcal{T})\cdot tr(\mathcal{F}(\mathcal{P}t+\mathcal{T}s))
-tr(\mathcal{F}\mathcal{T}\mathcal{F}(\mathcal{P}t+\mathcal{T}s)) \nonumber \\
&&+2tr(\mathcal{N}\mathcal{F}\mathcal{P}(\mathcal{P}t+\mathcal{T}s))+2tr(h\mathcal{T})\cdot tr(h(\mathcal{P}t+\mathcal{T}s)) \nonumber \\
&&-tr(h\mathcal{T}h(\mathcal{P}t+\mathcal{T}s))
+tr((h\mathcal{P}+s\mathcal{F})\mathcal{T})\cdot tr((h\mathcal{P}+s\mathcal{F})(\mathcal{P}t+\mathcal{T}s)) \nonumber \\
&&-tr((h\mathcal{P}+s\mathcal{F})\mathcal{T}(h\mathcal{P}+s\mathcal{F})(\mathcal{P}t+\mathcal{T}s)) \nonumber \\
&&+2tr((\mathcal{F}t+\mathcal{N}s)((h\mathcal{P}+s\mathcal{F})(f\mathcal{P}+t\mathcal{F})\mathcal{T})\Big\} \nonumber \\
&&+\frac{c_1-c_2}{16}\Big\{ (tr(\mathcal{F}(\mathcal{P}t+\mathcal{T}s)))^2
-tr(\mathcal{T}\mathcal{F}(\mathcal{P}t+\mathcal{T}s)(h\mathcal{P}+s\mathcal{F})) \nonumber \\
&&-tr(\mathcal{F}\mathcal{T}(h\mathcal{P}+s\mathcal{F})(\mathcal{P}t+\mathcal{T}s))
+tr(\mathcal{F}\mathcal{T})\cdot tr((h\mathcal{P}+s\mathcal{F})(\mathcal{P}t+\mathcal{T}s)) \nonumber \\
&&+2tr(\mathcal{T}\mathcal{N}(h\mathcal{P}+s\mathcal{F})(f\mathcal{P}+t\mathcal{F}))
+2tr((\mathcal{F}t+\mathcal{N}s)(h\mathcal{P}+s\mathcal{F})\mathcal{P}\mathcal{T})   \Big\}\nonumber
\end{eqnarray}
where
\begin{eqnarray}
&&-tr(\mathcal{T}\mathcal{F}(\mathcal{P}t+\mathcal{T}s)(h\mathcal{P}+s\mathcal{F})) \nonumber \\
&=&-\sum_\alpha g(\mathcal{F}(\mathcal{P}t+\mathcal{T}s)v_\alpha,\mathcal{F}(\mathcal{P}t+\mathcal{T}s)v_\alpha) \nonumber \\
&\leq & 0
\end{eqnarray}
and
\begin{eqnarray}
&&tr(\mathcal{F}\mathcal{T})\cdot tr((h\mathcal{P}+s\mathcal{F})(\mathcal{P}t+\mathcal{T}s)) \nonumber \\
&=&\sum_{\alpha,\beta}g(\mathcal{T}v_\alpha,\mathcal{T}v_\alpha)\cdot
g((\mathcal{P}t+\mathcal{T}s)v_\beta,(\mathcal{P}t+\mathcal{T}s)v_\beta) \nonumber \\
&\geq & 0
\end{eqnarray}

\begin{eqnarray}
&&g([A_{\mathcal{N}v_\alpha},A_\alpha]e_i,(f\mathcal{P}+t\mathcal{F})e_i)=2tr(A_\alpha A_{\mathcal{N}v_\alpha}(f\mathcal{P}+t\mathcal{F}))\\
&=&\frac{c_1+c_2}{16}\Big\{-2tr(\mathcal{T}\mathcal{N}\mathcal{F}(f\mathcal{P}+t\mathcal{F}))
-2tr(\mathcal{P}(f\mathcal{P}+t\mathcal{F}))\cdot tr(\mathcal{N}^2) \nonumber \\
&&+3tr(t\mathcal{N}h(f\mathcal{P}+t\mathcal{F}))-2tr((\mathcal{P}t+\mathcal{T}s)\mathcal{N}(h\mathcal{P}+s\mathcal{F})(f\mathcal{P}+t\mathcal{F})) \nonumber \\
&&-2tr((f\mathcal{P}+t\mathcal{F})^2)\cdot tr(\mathcal{N}(\mathcal{F}t+\mathcal{N}s))\Big\}   \nonumber\\
&&+\frac{c_1-c_2}{16}\Big\{-2tr(\mathcal{T}\mathcal{N}(h\mathcal{P}+s\mathcal{F})(f\mathcal{P}+t\mathcal{F}))
-2tr(\mathcal{N}\mathcal{F}(f\mathcal{P}+t\mathcal{F})(\mathcal{P}t+\mathcal{T}s))  \nonumber \\
&&-2tr((f\mathcal{P}+t\mathcal{F})^2)\cdot tr(\mathcal{N}^2)-2tr(\mathcal{P}(f\mathcal{P}+t\mathcal{F}))\cdot tr(\mathcal{N}(\mathcal{F}t+\mathcal{N}s)) \Big\}   \nonumber
\end{eqnarray}
where
\begin{eqnarray}
&&-2tr((f\mathcal{P}+t\mathcal{F})^2)\cdot tr(\mathcal{N}^2)  \nonumber\\
&=&-2g((f\mathcal{P}+t\mathcal{F})e_i,(f\mathcal{P}+t\mathcal{F})e_i)\cdot g(\mathcal{N}v_\alpha,\mathcal{N}v_\alpha) \nonumber\\
&\leq & 0
\end{eqnarray}

\section{Minimal submanifolds with flat normal connection}

In this section, we give the Simons' type integral formula for a compact minimal submanifold $M$
in $\overline{M}^m\times\overline{M}^n$. According to \cite{MK}, we have to compute the
following terms:
\begin{eqnarray}
&&\sum_i(\overline{\nabla}_X(\overline{R}(e_i,Y)e_i)^\bot)^\bot \nonumber\\
&=&\sum_i(\overline{R}(B(X,e_i),Y)e_i)^\bot+\sum_i(\overline{R}(e_i,B(X,Y))e_i)^\bot\\
&&+\sum_i(\overline{R}(e_i,Y)B(X,e_i))^\bot-\sum_iB(X,(\overline{R}(e_i,Y)e_i)^\top)\nonumber
\end{eqnarray}

where

\begin{eqnarray}
&&\sum_i(\overline{R}(B(X,e_i),Y)e_i)^\bot  \nonumber\\
&=&\frac{c_1+c_2}{16}\Big[B(X,Y)+\mathcal{N}B(X,\mathcal{P}Y)  \nonumber\\
&&+\sum_ig(B(X,e_i),\mathcal{F}e_i)\mathcal{F}Y+2\sum_ig(B(X,e_i),\mathcal{F}Y)\mathcal{F}e_i  \nonumber\\
&&+2sB(X,fY)-\sum_ig(B(X, e_i),he_i)hY +(\mathcal{F}t+\mathcal{N}s)B(X,(f\mathcal{P}+t\mathcal{F})Y) \nonumber\\
&&+\sum_ig(B(X,e_i),(h\mathcal{P}+s\mathcal{F})e_i)(h\mathcal{P}+s\mathcal{F})Y  \\
&&+2\sum_ig(B(X,e_i),(h\mathcal{P}+s\mathcal{F})Y)(h\mathcal{P}+s\mathcal{F})e_i\Big]  \nonumber\\
&&+\frac{c_1-c_2}{16}\Big[B(X,fY)+sB(X,Y)+\mathcal{N}B(X,(f\mathcal{P}+t\mathcal{F})Y)  \nonumber\\
&&+g(B(X,e_i),(h\mathcal{P}+s\mathcal{F})e_i)\mathcal{F}Y+(\mathcal{F}t+\mathcal{N}s)B(X,\mathcal{P}Y) \nonumber\\
&&+g(B(X,e_i),\mathcal{F}e_i)(h\mathcal{P}+s\mathcal{F})Y+2g(B(X,e_i),(h\mathcal{P}+s\mathcal{F})Y)\mathcal{F}e_i \nonumber\\
&&+2g(B(X,e_i),\mathcal{F}Y)(h\mathcal{P}+s\mathcal{F})e_i \Big]  \nonumber
\end{eqnarray}

\begin{eqnarray}
&&\sum_i(\overline{R}(e_i,B(X,Y))e_i)^\bot  \nonumber\\
&=&\frac{c_1+c_2}{16}\Big[-l\cdot B(X,Y)+3\mathcal{F}\mathcal{T}B(X,Y)+2htB(X,Y)\\
&&-tr(f)sB(X,Y)+3(h\mathcal{P}+s\mathcal{F})(\mathcal{P}t+\mathcal{T}s)B(X,Y)\Big] \nonumber\\
&&+\frac{c_1-c_2}{16}\Big[-tr(f)B(X,Y)-l\cdot sB(X,Y)\nonumber\\
&&+3\mathcal{F}(\mathcal{P}t+\mathcal{T}s)B(X,Y)+3(h\mathcal{P}+s\mathcal{F})\mathcal{T}B(X,Y)     \Big]\nonumber
\end{eqnarray}

\begin{eqnarray}
&&\sum_i(\overline{R}(e_i,Y)B(X,e_i))^\bot  \nonumber\\
&=&\frac{c_1+c_2}{16}\Big[\sum_ig(\mathcal{F}Y,B(X,e_i))\mathcal{F}e_i-\sum_ig(\mathcal{F}e_i,B(X,e_i))\mathcal{F}Y  \nonumber\\
&&+2\mathcal{N}B(X,\mathcal{P}Y)+2\sum_ig(hY,B(X,e_i))he_i-\sum_ig(he_i,B(X,e_i))hY   \nonumber\\
&&+\sum_ig((h\mathcal{P}+s\mathcal{F})Y,B(X,e_i))(h\mathcal{P}+s\mathcal{F})e_i    \\
&&-\sum_ig((h\mathcal{P}+s\mathcal{F})e_i,B(X,e_i))(h\mathcal{P}+s\mathcal{F})Y  \nonumber\\
&&+2(\mathcal{F}t+\mathcal{N}s)B(X,(f\mathcal{P}+t\mathcal{F})Y)  \Big]\nonumber\\
&&+\frac{c_1-c_2}{16}\Big[g((h\mathcal{P}+s\mathcal{F})Y,B(X,e_i))\mathcal{F}e_i-g((h\mathcal{P}+s\mathcal{F})e_i,B(X,e_i))\mathcal{F}Y  \nonumber\\
&&+g(\mathcal{F}Y,B(X,e_i))(h\mathcal{P}+s\mathcal{F})e_i-g(\mathcal{F}e_i,B(X,e_i))(h\mathcal{P}+s\mathcal{F})Y \nonumber\\
&&+2\mathcal{N}B(X,(f\mathcal{P}+t\mathcal{F})Y)+2(\mathcal{F}t+\mathcal{N}s)B(X,\mathcal{P}Y) \Big]\nonumber
\end{eqnarray}

\begin{eqnarray}
&&\sum_i(\overline{R}(e_i,Y)e_i)^\top  \\
&=&\frac{c_1+c_2}{16}\Big[(1-l)Y+3\mathcal{P}^2Y+2f^2Y-tr(f)\cdot fY+3(f\mathcal{P}+t\mathcal{F})^2Y \Big]\nonumber\\
&&+\frac{c_1-c_2}{16}\Big[(2-l)fY-tr(f)Y+3\mathcal{P}(f\mathcal{P}+t\mathcal{F})Y+3(f\mathcal{P}+t\mathcal{F})\mathcal{P}Y \Big]\nonumber
\end{eqnarray}
Therefore
\begin{eqnarray}
&&-\sum_iB(X,(\overline{R}(e_i,Y)e_i)^\top)  \\
&=&-\frac{c_1+c_2}{16}\Big[(1-l)B(X,Y)+3B(X,\mathcal{P}^2Y)+2B(X,f^2Y) \nonumber\\
&&-tr(f)B(X,fY)+3B(X,(f\mathcal{P}+t\mathcal{F})^2Y) \Big]  \nonumber\\
&&-\frac{c_1-c_2}{16}\Big[(2-l)B(X,fY)-tr(f)B(X,Y)+3B(X,\mathcal{P}(f\mathcal{P}+t\mathcal{F})Y)  \nonumber\\
&&+3B(X,(f\mathcal{P}+t\mathcal{F})\mathcal{P}Y) \Big]\nonumber
\end{eqnarray}

We also have to compute
\begin{eqnarray}
&&\sum_i(\overline{\nabla}_{e_i}(\overline{R}(e_i,X)Y)^\bot)^\bot  \\
&=&\sum_i(\overline{R}(e_i,B(e_i,X))Y)^\bot
+\sum_i(\overline{R}(e_i,X)B(e_i,Y))^\bot-\sum_iB(e_i,(\overline{R}(e_i,X)Y)^\top)\nonumber
\end{eqnarray}

where
\begin{eqnarray}
&&\sum_i(\overline{R}(e_i,B(e_i,X))Y)^\bot  \nonumber\\
&=&\frac{c_1+c_2}{16}\Big[-B(X,Y)-sB(X,fY)+\mathcal{N}B(X,\mathcal{P}Y)\nonumber\\
&&-\sum_ig(B(e_i,X),\mathcal{F}Y)\mathcal{F}e_i-2\sum_ig(B(e_i,X),\mathcal{F}e_i)\mathcal{F}Y  \nonumber\\
&&+2\sum_ig(B(e_i,X),hY)he_i+(\mathcal{F}t+\mathcal{N}s)B((f\mathcal{P}+t\mathcal{F})Y,X) \nonumber \\
&&-\sum_ig(B(e_i,X),(h\mathcal{P}+s\mathcal{F})Y)(h\mathcal{P}+s\mathcal{F})e_i  \\
&&-2\sum_ig(B(e_i,X),(h\mathcal{P}+s\mathcal{F})e_i)(h\mathcal{P}+s\mathcal{F})Y \Big] \nonumber\\
&&+\frac{c_1-c_2}{16}\Big[-B(X,fY)-hB(X,Y)+\mathcal{N}B(X,(f\mathcal{P}+t\mathcal{F})Y)  \nonumber\\
&&+(\mathcal{F}t+\mathcal{N}s)B(X,\mathcal{P}Y)   \nonumber\\
&&-g(B(e_i,X),(h\mathcal{P}+s\mathcal{F})Y)\mathcal{F}e_i-2g(B(e_i,X),(h\mathcal{P}+s\mathcal{F})e_i)\mathcal{F}Y  \nonumber\\
&&-g(B(e_i,X),\mathcal{F}Y)(h\mathcal{P}+s\mathcal{F})e_i-2g(B(e_i,X),\mathcal{F}e_i)(h\mathcal{P}+s\mathcal{F})Y \Big] \nonumber
\end{eqnarray}

\begin{eqnarray}
&&\sum_i(\overline{R}(e_i,X)B(Y,e_i))^\bot  \nonumber\\
&=&\frac{c_1+c_2}{16}\Big[\sum_ig(\mathcal{F}X,B(Y,e_i))\mathcal{F}e_i-\sum_ig(\mathcal{F}e_i,B(Y,e_i))\mathcal{F}X  \nonumber\\
&&+2\mathcal{N}B(Y,\mathcal{P}X)+2\sum_ig(hX,B(Y,e_i))he_i-\sum_ig(he_i,B(Y,e_i))hX   \nonumber\\
&&+\sum_ig((h\mathcal{P}+s\mathcal{F})X,B(Y,e_i))(h\mathcal{P}+s\mathcal{F})e_i    \\
&&-\sum_ig((h\mathcal{P}+s\mathcal{F})e_i,B(Y,e_i))(h\mathcal{P}+s\mathcal{F})X \nonumber\\
&&+2(\mathcal{F}t+\mathcal{N}s)B(Y,(f\mathcal{P}+t\mathcal{F})X)  \Big]\nonumber\\
&&+\frac{c_1-c_2}{16}\Big[g((h\mathcal{P}+s\mathcal{F})X,B(Y,e_i))\mathcal{F}e_i-g((h\mathcal{P}+s\mathcal{F})e_i,B(Y,e_i))\mathcal{F}X \nonumber\\
&&+g(\mathcal{F}X,B(Y,e_i))(h\mathcal{P}+s\mathcal{F})e_i-g(\mathcal{F}e_i,B(Y,e_i))(h\mathcal{P}+s\mathcal{F})X \nonumber\\
&&+2\mathcal{N}B(Y,(f\mathcal{P}+t\mathcal{F})X)+2(\mathcal{F}t+\mathcal{N}s)B(Y,\mathcal{P}X)   \Big]\nonumber
\end{eqnarray}

\begin{eqnarray}
&&\sum_i(\overline{R}(e_i,X)Y)^\top \nonumber\\
&=&\frac{c_1+c_2}{16}\sum_i\Big[g(X,Y)e_i-g(e_i,Y)X+g(\mathcal{P}X,Y)\mathcal{P}e_i+g(\mathcal{P}Y,e_i)\mathcal{P}X  \nonumber\\
&&+2g(e_i,\mathcal{P}X)\mathcal{P}Y+2g(fX,Y)fe_i-g(e_i,fY)fX   \\
&&+g((f\mathcal{P}+t\mathcal{F})X,Y)(f\mathcal{P}+t\mathcal{F})e_i  \nonumber\\
&&+g((f\mathcal{P}+t\mathcal{F})Y,e_i)(f\mathcal{P}+t\mathcal{F})X
+2g((f\mathcal{P}+t\mathcal{F})X,e_i)(f\mathcal{P}+t\mathcal{F})Y \Big]  \nonumber\\
&&+\frac{c_1-c_2}{16}\sum_i\Big[g(fX,Y)e_i-g(fY,e_i)X+g(X,Y)fe_i-g(Y,e_i)fX \nonumber\\
&&+g((f\mathcal{P}+t\mathcal{F})X,Y)\mathcal{P}e_i+g((f\mathcal{P}+t\mathcal{F})Y,e_i)\mathcal{P}X \nonumber\\
&&+g(\mathcal{P}X,Y)(f\mathcal{P}+t\mathcal{F})e_i+g(\mathcal{P}Y,e_i)(f\mathcal{P}+t\mathcal{F})X \nonumber\\
&&+2g((f\mathcal{P}+t\mathcal{F})X,e_i)\mathcal{P}Y+2g(\mathcal{P}X,e_i)(f\mathcal{P}+t\mathcal{F})Y  \Big]  \nonumber
\end{eqnarray}

Therefore
\begin{eqnarray}
&&-\sum_iB(e_i,(\overline{R}(e_i,X)Y)^\top)   \\
&&=-\frac{c_1+c_2}{16}\Big[-B(X,Y)+3B(\mathcal{P}X,\mathcal{P}Y)+2g(fX,Y)B(fe_i,e_i) \nonumber\\
&&-B(fX,fY)+3B((f\mathcal{P}+t\mathcal{F})X,(f\mathcal{P}+t\mathcal{F})Y)\Big] \nonumber\\
&&-\frac{c_1-c_2}{16}\Big[-B(X,fY)+g(X,Y)B(e_i,fe_i)-B(fX,Y)   \nonumber\\
&&+3B(\mathcal{P}X,(f\mathcal{P}+t\mathcal{F})Y)+3B(\mathcal{P}Y,(f\mathcal{P}+t\mathcal{F})X) \Big]  \nonumber
\end{eqnarray}

Then we have\\

\textbf{Lemma 3.1} Let $M$ be a minimal submanifold in $\overline{M}^m\times\overline{M}^n$ with flat normal connection,
then
\begin{eqnarray}
&&g(\nabla^2B,B)=\sum_{i,j,k}g(\nabla_{e_i}\nabla_{e_i}B(e_j,e_k),B(e_j,e_k))  \nonumber\\
&=&\sum_{i,j,k}\Big[g((R(e_i,e_j)B)(e_i,e_k),B(e_j,e_k)) \nonumber\\
&&+g((\overline{\nabla}_{e_j}(\overline{R}(e_i,e_k)e_i)^\bot)^\bot,B(e_j,e_k))
+g((\overline{\nabla}_{e_i}(\overline{R}(e_i,e_j)e_k)^\bot)^\bot,B(e_j,e_k))\nonumber \\
&=&\sum_{i,j,k}g((R(e_i,e_j)B)(e_i,e_k),B(e_j,e_k)) \\
&&+\frac{c_1+c_2}{16}\Big\{3\sum_{\alpha}\Big[tr(A_\alpha A_{\mathcal{F}\mathcal{T}v_\alpha})
-2tr(A_\alpha A_{\mathcal{N}v_\alpha}\mathcal{P})
+tr((A_\alpha\mathcal{P})^2)-tr(A_\alpha^2\mathcal{P}^2)\Big] \nonumber\\
&&+3\sum_{\alpha,\beta}\Big[g(A_\alpha(\mathcal{P}t+\mathcal{T}s)v_\beta, A_\beta(\mathcal{P}t+\mathcal{T}s)v_\alpha)
-g(A_\alpha(\mathcal{P}t+\mathcal{T}s)v_\alpha, A_\beta(\mathcal{P}t+\mathcal{T}s)v_\beta)  \Big] \nonumber\\
&&+3\sum_{\alpha,\beta}\Big[g(A_\alpha\mathcal{T}v_\beta, A_\beta\mathcal{T}v_\alpha)
-g(A_\alpha\mathcal{T}v_\alpha, A_\beta\mathcal{T}v_\beta)\Big] \nonumber\\
&&+3\sum_{\alpha}\Big[tr((A_\alpha(f\mathcal{P}+t\mathcal{F}))^2)-tr(A_\alpha^2(f\mathcal{P}+t\mathcal{F})^2)\Big] \nonumber\\
&&+3\sum_{\alpha,\beta}\Big[g(A_\alpha tv_\beta, A_\beta tv_\alpha)-g(A_\alpha tv_\alpha, A_\beta tv_\beta)\Big]
+3\sum_{\alpha,\beta}g(A_\alpha tv_\beta, A_\beta tv_\alpha)  \nonumber\\
&&+\sum_{\alpha}\Big[2tr(A_\alpha A_{htv_\alpha})+tr(A_\alpha A_{sv_\alpha}f)-tr(f)tr(A_\alpha A_{sv_\alpha})+tr(f)tr(A_\alpha^2f) \nonumber\\
&&-2tr(A_\alpha^2f^2)-2(tr(A_\alpha f))^2+tr((A_\alpha f)^2) \nonumber\\
&&-6tr(A_\alpha A_{(\mathcal{F}t+\mathcal{N}s)v_\alpha}(f\mathcal{P}+t\mathcal{F}))
+3tr(A_\alpha A_{(h\mathcal{P}+s\mathcal{F})(\mathcal{P}t+\mathcal{T}s)v_\alpha})\Big] \Big\}\nonumber
\end{eqnarray}
\begin{eqnarray}
&&+\frac{c_1-c_2}{16}\Big\{l\cdot tr(A_\alpha^2f)-l\cdot tr(A_\alpha A_{sv_\alpha}) \nonumber\\
&&+6g(A_\alpha\mathcal{T}v_\beta,A_\beta(\mathcal{P}t+\mathcal{T}s)v_\alpha)
-6g(A_\alpha\mathcal{T}v_\alpha,A_\beta(\mathcal{P}t+\mathcal{T}s)v_\beta)\nonumber\\
&&-6tr(A_\alpha A_{\mathcal{N}v_\alpha}(f\mathcal{P}+t\mathcal{F}))-6tr(A_\alpha A_{(\mathcal{F}t+\mathcal{N}s)}v_\alpha\mathcal{P}) \nonumber\\              &&+6tr(A_\alpha A_{(h\mathcal{P}+s\mathcal{F})\mathcal{T}v_\alpha})-6tr(A_\alpha^2(f\mathcal{P}+t\mathcal{F})\mathcal{P})
+6tr(A_\alpha\mathcal{P}A_\alpha(f\mathcal{P}+t\mathcal{F}))     \Big\}\nonumber
\end{eqnarray}
$\Box$\vspace{2mm}

We also have the following relations:
\begin{eqnarray}
\sum_{\alpha}tr(A_\alpha A_{\mathcal{F}\mathcal{T}v_\alpha})=\sum_{\alpha}[tr(A_{\mathcal{N}v_\alpha}^2)-tr(A_\alpha^2)]
\end{eqnarray}
\begin{eqnarray}
\sum_{\alpha}tr(A_\alpha A_{htv_\alpha})&=&\sum_{\alpha}[tr(A_{\alpha}^2)-tr(A_{s^2v_\alpha})]\nonumber\\
&=&\sum_{\alpha}[tr(A_{\alpha}^2)-tr(A_{sv_\alpha}^2)]
\end{eqnarray}
\begin{eqnarray}
\sum_{\alpha} tr(A_{s^2v_\alpha})=\sum_{\alpha} tr(A_{sv_\alpha}^2)
\end{eqnarray}
\begin{eqnarray}
\sum_{\alpha}|[\mathcal{P}, A_\alpha]|^2=2\sum_{\alpha}[tr((A_\alpha\mathcal{P})^2)-tr(A_\alpha^2\mathcal{P}^2)]
\end{eqnarray}
\begin{eqnarray}
\sum_{\alpha}|[f, A_\alpha]|^2=-2\sum_{\alpha}[tr((A_\alpha f)^2)-tr(A_\alpha^2 f^2)]
\end{eqnarray}
\begin{eqnarray}
&&\sum_{\alpha}|[f\mathcal{P}+t\mathcal{F}, A_\alpha]|^2 \nonumber\\
&=&2\sum_{\alpha}[tr((A_\alpha (f\mathcal{P}+t\mathcal{F}))^2)-tr(A_\alpha^2 (f\mathcal{P}+t\mathcal{F})^2)]
\end{eqnarray}
\begin{eqnarray}
\sum_{\alpha}tr(A_\alpha A_{\mathcal{N}\alpha})=\sum_{\alpha}tr(A_\alpha^2\mathcal{P})=\sum_{\alpha}tr(A_\alpha\mathcal{P})=0
\end{eqnarray}
\begin{eqnarray}
\sum_{\alpha}tr(A_\alpha(f\mathcal{P}+t\mathcal{F}))=\sum_{\alpha}tr(A_\alpha^2(f\mathcal{P}+t\mathcal{F}))=0
\end{eqnarray}
\begin{eqnarray}
\sum_{\alpha}tr(A_\alpha A_{(\mathcal{F}t+\mathcal{N}s)v_\alpha})=0
\end{eqnarray}
\begin{eqnarray}
&&\sum_{\alpha}tr(A_\alpha A_{(h\mathcal{P}+s\mathcal{F})(\mathcal{P}t+\mathcal{T}s)v_\alpha}) \nonumber\\
&=&-\sum_{\alpha,\beta}tr(A_\alpha A_\beta)\cdot g((\mathcal{P}t+\mathcal{T}s)v_{\alpha},(\mathcal{P}t+\mathcal{T}s)v_{\beta}) \nonumber\\
&=&-\sum_{i,j,k}\Big[g((h\mathcal{P}+s\mathcal{F})e_j,B(e_i,e_k))\Big]^2\leq 0
\end{eqnarray}
\begin{eqnarray}
tr(A_\alpha^2\mathcal{P}(f\mathcal{P}+t\mathcal{F}))=tr(A_\alpha^2(f\mathcal{P}+t\mathcal{F})\mathcal{P})
\end{eqnarray}

Now we can rewrite (3.12) as follows:
\begin{eqnarray}
&&g(\nabla^2A,A)=g(\nabla^2B,B)=\sum_{i,j,k}g(\nabla_{e_i}\nabla_{e_i}B(e_j,e_k),B(e_j,e_k))  \nonumber\\
&=&\sum_{i,j,\alpha}g((R(e_i,e_j)A_\alpha)e_i,A_\alpha e_j) \\
&&+\frac{c_1+c_2}{16}\Big\{3\sum_{\alpha}\Big[tr(A_{\mathcal{N}v_\alpha}^2)-tr(A_\alpha^2)
-2tr(A_\alpha A_{\mathcal{N}v_\alpha}\mathcal{P})
+\frac{1}{2}\sum_\alpha|[\mathcal{P},A_\alpha]|^2\Big] \nonumber\\
&&+3\frac{c_1+c_2}{16}\Big[(tr(\mathcal{F}(\mathcal{P}t+\mathcal{T}s)))^2-tr((\mathcal{F}(\mathcal{P}t+\mathcal{T}s))^2)  \nonumber\\
&&+2tr(\mathcal{N}(h\mathcal{P}+s\mathcal{F})\mathcal{P}(\mathcal{P}t+\mathcal{T}s))+2(tr(h(\mathcal{P}t+\mathcal{T}s)))^2
-tr((h(\mathcal{P}t+\mathcal{T}s))^2)  \nonumber\\
&&+(tr((h\mathcal{P}+s\mathcal{F})(\mathcal{P}t+\mathcal{T}s)))^2-tr(((h\mathcal{P}+s\mathcal{F})(\mathcal{P}t+\mathcal{T}s))^2)  \nonumber\\
&&+2tr((\mathcal{F}t+\mathcal{N}s)(h\mathcal{P}+s\mathcal{F})(f\mathcal{P}+t\mathcal{F})(\mathcal{P}t+\mathcal{T}s)) \Big]\nonumber\\
&&+6\frac{c_1-c_2}{16}\Big[tr(\mathcal{F}(\mathcal{P}t+\mathcal{T}s))\cdot tr((h\mathcal{P}+s\mathcal{F})(\mathcal{P}t+\mathcal{T}s))\nonumber\\
&&-tr(\mathcal{F}(\mathcal{P}t+\mathcal{T}s)(h\mathcal{P}+s\mathcal{F})(\mathcal{P}t+\mathcal{T}s))
+tr(\mathcal{N}(h\mathcal{P}+s\mathcal{F})(f\mathcal{P}+t\mathcal{F})(\mathcal{P}t+\mathcal{T}s))  \nonumber\\
&&+tr((\mathcal{F}t+\mathcal{N}s)(h\mathcal{P}+s\mathcal{F})\mathcal{P}(\mathcal{P}t+\mathcal{T}s))  \Big]\nonumber\\
&&+3\frac{c_1+c_2}{16}\Big[(tr(\mathcal{F}\mathcal{T}))^2-tr((\mathcal{F}\mathcal{T})^2)
+2tr(\mathcal{N}\mathcal{F}\mathcal{P}\mathcal{T})-tr((\mathcal{F}t)^2)  \nonumber\\
&&+(tr(h\mathcal{P}+s\mathcal{F})\mathcal{T})^2-tr(((h\mathcal{P}+s\mathcal{F})\mathcal{T})^2)
+2tr((\mathcal{F}t+\mathcal{N}s)\mathcal{F}(f\mathcal{P}+t\mathcal{F})\mathcal{T}) \Big]  \nonumber\\
&&+6\frac{c_1-c_2}{16}\Big[tr(\mathcal{F}\mathcal{T})tr((h\mathcal{P}+s\mathcal{F})\mathcal{T})
-tr(\mathcal{F}\mathcal{T}(h\mathcal{P}+s\mathcal{F})\mathcal{T})  \nonumber\\
&&+tr(\mathcal{N}\mathcal{F}(f\mathcal{P}+t\mathcal{F})\mathcal{T})+tr(\mathcal{F}\mathcal{P}\mathcal{T}(\mathcal{F}t+\mathcal{N}s)) \Big]  \nonumber\\
&&+3\frac{c_1+c_2}{16}\Big[-tr((\mathcal{F}t)^2)-2tr(\mathcal{N}h\mathcal{P}t)+2(tr(ht))^2-tr((ht)^2)\nonumber
\end{eqnarray}
\begin{eqnarray}
&&+(tr((h\mathcal{P}+s\mathcal{F})t))^2-tr(((h\mathcal{P}+s\mathcal{F})t)^2)-2tr((\mathcal{F}t+\mathcal{N}s)h(f\mathcal{P}+t\mathcal{F})t)\Big]\nonumber\\
&&+6\frac{c_1-c_2}{16}\Big[-tr((h\mathcal{P}+s\mathcal{F})t\mathcal{F}t)
-tr((f\mathcal{P}+t\mathcal{F})t\mathcal{N}h)-tr(h\mathcal{P}t(\mathcal{F}t+\mathcal{N}s))  \Big]  \nonumber\\
&&-3\frac{c_1+c_2}{16}\Big[-2tr((f\mathcal{P}+t\mathcal{F})\mathcal{T}(\mathcal{F}t+\mathcal{N}s)\mathcal{F})
-2tr(\mathcal{P}(f\mathcal{P}+t\mathcal{F}))tr(\mathcal{N}(\mathcal{F}t+\mathcal{N}s))  \nonumber\\
&&+3tr(t(\mathcal{F}t+\mathcal{N}s)h(f\mathcal{P}+t\mathcal{F}))
-2tr((\mathcal{P}t+\mathcal{T}s)(\mathcal{F}t+\mathcal{N}s)(h\mathcal{P}+s\mathcal{F})(f\mathcal{P}+t\mathcal{F})) \nonumber\\
&&-2tr((f\mathcal{P}+t\mathcal{F})^2)\cdot tr((\mathcal{F}t+\mathcal{N}s)^2) \Big] \nonumber\\
&&+6\frac{c_1-c_2}{16}\Big[tr(\mathcal{T}(\mathcal{F}t+\mathcal{N}s)(h\mathcal{P}+s\mathcal{F})(f\mathcal{P}+t\mathcal{F})) \nonumber\\
&&+tr((\mathcal{P}t+\mathcal{T}s)(\mathcal{F}t+\mathcal{N}s)\mathcal{F}(f\mathcal{P}+t\mathcal{F})) \nonumber\\
&&+tr((f\mathcal{P}+t\mathcal{F})^2)\cdot tr(\mathcal{N}(\mathcal{F}t+\mathcal{N}s))
+tr((\mathcal{F}t+\mathcal{N}s)^2)\cdot tr(\mathcal{P}(f\mathcal{P}+t\mathcal{F})) \Big]   \nonumber \\
&&+\frac{3}{2}\sum_{\alpha}|[f\mathcal{P}+t\mathcal{F},A_\alpha]|^2
+3\sum_{\alpha,\beta}g(A_\alpha tv_\beta, A_\beta tv_\alpha)  \nonumber\\
&&+\sum_{\alpha}\Big[2tr(A_\alpha^2)-2tr(A_{sv_{\alpha}}^2)
+tr(A_\alpha A_{sv_\alpha}f)-tr(f)tr(A_\alpha A_{sv_\alpha})+tr(f)tr(A_\alpha^2f) \nonumber\\
&&-2tr(A_\alpha^2f^2)-2(tr(A_\alpha f))^2+tr((A_\alpha f)^2)
+3tr(A_\alpha A_{(h\mathcal{P}+s\mathcal{F})(\mathcal{P}t+\mathcal{T}s)v_\alpha})\Big] \Big\}\nonumber\\
&&+\frac{c_1-c_2}{16}\Big\{l\cdot tr(A_\alpha^2f)-l\cdot tr(A_\alpha A_{sv_\alpha}) \nonumber\\
&&+6\frac{c_1+c_2}{16}\Big[tr(\mathcal{F}\mathcal{T})\cdot tr(\mathcal{F}(\mathcal{P}t+\mathcal{T}s))
-tr(\mathcal{F}\mathcal{T}\mathcal{F}(\mathcal{P}t+\mathcal{T}s)) \nonumber \\
&&+2tr(\mathcal{N}\mathcal{F}\mathcal{P}(\mathcal{P}t+\mathcal{T}s))+2tr(h\mathcal{T})\cdot tr(h(\mathcal{P}t+\mathcal{T}s)) \nonumber \\
&&-tr(h\mathcal{T}h(\mathcal{P}t+\mathcal{T}s))
+tr((h\mathcal{P}+s\mathcal{F})\mathcal{T})\cdot tr((h\mathcal{P}+s\mathcal{F})(\mathcal{P}t+\mathcal{T}s)) \nonumber \\
&&-tr((h\mathcal{P}+s\mathcal{F})\mathcal{T}(h\mathcal{P}+s\mathcal{F})(\mathcal{P}t+\mathcal{T}s)) \nonumber \\
&&+2tr((\mathcal{F}t+\mathcal{N}s)((h\mathcal{P}+s\mathcal{F})(f\mathcal{P}+t\mathcal{F})\mathcal{T})\Big] \nonumber \\
&&+6\frac{c_1-c_2}{16}\Big[ (tr(\mathcal{F}(\mathcal{P}t+\mathcal{T}s)))^2
-tr(\mathcal{T}\mathcal{F}(\mathcal{P}t+\mathcal{T}s)(h\mathcal{P}+s\mathcal{F})) \nonumber \\
&&-tr(\mathcal{F}\mathcal{T}(h\mathcal{P}+s\mathcal{F})(\mathcal{P}t+\mathcal{T}s))
+tr(\mathcal{F}\mathcal{T})\cdot tr((h\mathcal{P}+s\mathcal{F})(\mathcal{P}t+\mathcal{T}s)) \nonumber \\
&&+2tr(\mathcal{T}\mathcal{N}(h\mathcal{P}+s\mathcal{F})(f\mathcal{P}+t\mathcal{F}))
+2tr((\mathcal{F}t+\mathcal{N}s)(h\mathcal{P}+s\mathcal{F})\mathcal{P}\mathcal{T}) \Big] \nonumber\\
&&-3\frac{c_1+c_2}{16}\Big[-2tr(\mathcal{T}\mathcal{N}\mathcal{F}(f\mathcal{P}+t\mathcal{F}))
-2tr(\mathcal{P}(f\mathcal{P}+t\mathcal{F}))\cdot tr(\mathcal{N}^2) \nonumber \\
&&+3tr(t\mathcal{N}h(f\mathcal{P}+t\mathcal{F}))-2tr((\mathcal{P}t+\mathcal{T}s)\mathcal{N}(h\mathcal{P}+s\mathcal{F})(f\mathcal{P}+t\mathcal{F})) \nonumber \\
&&-2tr((f\mathcal{P}+t\mathcal{F})^2)\cdot tr(\mathcal{N}(\mathcal{F}t+\mathcal{N}s))\Big]   \nonumber\\
&&-3\frac{c_1-c_2}{16}\Big[-2tr(\mathcal{T}\mathcal{N}(h\mathcal{P}+s\mathcal{F})(f\mathcal{P}+t\mathcal{F}))
-2tr(\mathcal{N}\mathcal{F}(f\mathcal{P}+t\mathcal{F})(\mathcal{P}t+\mathcal{T}s))  \nonumber \\
&&-2tr((f\mathcal{P}+t\mathcal{F})^2)\cdot tr(\mathcal{N}^2)
-2tr(\mathcal{P}(f\mathcal{P}+t\mathcal{F}))\cdot tr(\mathcal{N}(\mathcal{F}t+\mathcal{N}s)) \Big] \nonumber
\end{eqnarray}
\begin{eqnarray}
&&-3\frac{c_1+c_2}{16}\Big[-2tr(\mathcal{F}\mathcal{P}\mathcal{T}(\mathcal{F}t+\mathcal{N}s))
-2tr(\mathcal{P}^2)tr(\mathcal{N}(\mathcal{F}t+\mathcal{N}s))  \nonumber\\
&&+3tr(h\mathcal{P}t(\mathcal{F}t+\mathcal{N}s))
-2tr((\mathcal{P}t+\mathcal{T}s)(\mathcal{F}t+\mathcal{N}s)(h\mathcal{P}+s\mathcal{F})\mathcal{P})  \nonumber\\
&&-2tr((f\mathcal{P}+t\mathcal{F})\mathcal{P})\cdot tr((\mathcal{F}t+\mathcal{N}s)^2)\Big] \nonumber\\
&&-3\frac{c_1-c_2}{16}\Big[-4tr(\mathcal{P}\mathcal{T}(\mathcal{F}t+\mathcal{N}s)(h\mathcal{P}+s\mathcal{F}))  \nonumber\\
&&-2tr((f\mathcal{P}+t\mathcal{F})\mathcal{P})\cdot tr(\mathcal{N}(\mathcal{F}t+\mathcal{N}s))
-2tr(\mathcal{P}^2)\cdot tr((\mathcal{F}t+\mathcal{N}s)^2) \Big] \nonumber\\
&&+6tr(A_\alpha A_{(h\mathcal{P}+s\mathcal{F})\mathcal{T}v_\alpha})-6tr(A_\alpha^2(f\mathcal{P}+t\mathcal{F})\mathcal{P})
+6tr(A_\alpha\mathcal{P}A_\alpha(f\mathcal{P}+t\mathcal{F}))     \Big\}\nonumber
\end{eqnarray}

As in \cite{MK}, we have the following result:

\textbf{Lemma 3.2} Let $M$ be an $l$-dimensional minimal submanifold in $\overline{M}^m\times \overline{M}^n$.
If $U$ is a parallel section in the normal bundle of $M$, then
\begin{eqnarray}
&&div(\nabla_{\mathcal{T}U}\mathcal{T}U) \nonumber\\
&=&\frac{c_1+c_2}{16}\Big\{(l-1)g(\mathcal{T}U,\mathcal{T}U)
+3g(\mathcal{P}\mathcal{T}U,\mathcal{P}\mathcal{T}U)+2tr(f)g(f\mathcal{T}U,\mathcal{T}U)  \nonumber\\
&&-g(f\mathcal{T}U,f\mathcal{T}U)+3g((f\mathcal{P}+t\mathcal{F})\mathcal{T}U,(f\mathcal{P}+t\mathcal{F})\mathcal{T}U)\Big\}   \nonumber\\
&&+\frac{c_1-c_2}{16}\Big\{(l-2)g(f\mathcal{T}U,\mathcal{T}U)+tr(f)\cdot g(\mathcal{T}U,\mathcal{T}U) \nonumber\\
&&+6g(\mathcal{P}\mathcal{T}U,(f\mathcal{P}+t\mathcal{F})\mathcal{T}U)    \Big\} \nonumber\\
&&+\frac{1}{2}|[\mathcal{P},A_U]|^2+tr(A_{\mathcal{N}U}^2)-2tr(A_U A_{\mathcal{N}U}\mathcal{P})-tr(A_U^2)  \nonumber\\
&&+\sum_\alpha g(A_U\mathcal{T}v_\alpha,A_U\mathcal{T}v_\alpha)-\sum_\alpha g(A_\alpha\mathcal{T}U,A_\alpha\mathcal{T}U)
\end{eqnarray}

\textbf{Proof} From Gauss equation, we get
\begin{eqnarray}
&&S(\mathcal{T}U,\mathcal{T}U)  \nonumber\\
&=&\frac{c_1+c_2}{16}\Big\{(l-1)g(\mathcal{T}U,\mathcal{T}U)
+3g(\mathcal{P}\mathcal{T}U,\mathcal{P}\mathcal{T}U)+2tr(f)g(f\mathcal{T}U,\mathcal{T}U)  \nonumber\\
&&-g(f\mathcal{T}U,f\mathcal{T}U)+3g((f\mathcal{P}+t\mathcal{F})\mathcal{T}U,(f\mathcal{P}+t\mathcal{F})\mathcal{T}U)\Big\}  \nonumber\\
&&+\frac{c_1-c_2}{16}\Big\{(l-2)g(f\mathcal{T}U,\mathcal{T}U)+tr(f)\cdot g(\mathcal{T}U,\mathcal{T}U) \nonumber\\
&&+6g(\mathcal{P}\mathcal{T}U,(f\mathcal{P}+t\mathcal{F})\mathcal{T}U)    \Big\}
-\sum_\alpha g(A_\alpha\mathcal{T}U,A_\alpha\mathcal{T}U)
\end{eqnarray}
where $S$ is the Ricci tensor of $M$.

From \cite{MK} we know that
\begin{eqnarray}
\nabla_X(\mathcal{T}U)=-\mathcal{P}A_UX+A_{\mathcal{N}U}X
\end{eqnarray}
Therefore
\begin{eqnarray}
div(\mathcal{T}U)&=&\sum_i g(e_i,\nabla_{e_i}(\mathcal{T}U))=\sum_i g(e_i,-\mathcal{P}A_U e_i+A_{\mathcal{N}U}e_i)  \nonumber\\
&=&-tr(\mathcal{P}A_U)+tr(A_{\mathcal{N}U}) =0,
\end{eqnarray}
\begin{eqnarray}
&&|\nabla\mathcal{T}U|^2=\sum_i g(\nabla_{e_i}\mathcal{T}U,\nabla_{e_i}\mathcal{T}U)  \\
&=&tr(A_U^2)+tr(A_{\mathcal{N}U}^2)-2tr(A_U A_{\mathcal{N}U}\mathcal{P})-\sum_\alpha g(A_U \mathcal{T}v_\alpha,A_U \mathcal{T}v_\alpha),\nonumber
\end{eqnarray}
\begin{eqnarray}
&&|L_{\mathcal{T}U}g|^2=\sum_{i,j}\Big[L_{\mathcal{T}U}g(e_i,e_j)\Big]^2=\sum_{i,j}\Big[g(\nabla_{e_i}\mathcal{T}U,e_j)
+g(\nabla_{e_j}\mathcal{T}U,e_i) \Big]^2  \nonumber\\
&=&|[\mathcal{P},A_U]|^2+4tr(A_{\mathcal{N}U}^2)-8tr(A_U A_{\mathcal{N}U}\mathcal{P})
\end{eqnarray}
Put (3.26)-(3.30) into
\begin{eqnarray}
div(\nabla_{\mathcal{T}U}\mathcal{T}U)&=&div(div(\mathcal{T}U)\mathcal{T}U)+S(\mathcal{T}U,\mathcal{T}U)  \nonumber\\
&&+\frac{1}{2}|L_{\mathcal{T}U}g|^2-|\nabla\mathcal{T}U|^2-(div(\mathcal{T}U))^2 \nonumber
\end{eqnarray}
we get the result. $\Box$\vspace{2mm}

Since the normal connection of $M$ is flat, we can choose an orthonormal basis $v_\alpha$ of $T(M)^\bot$ such that
$Dv_\alpha=0$ for all $\alpha$. Thus
\begin{eqnarray}
&&div(\nabla_{\mathcal{T}v_\alpha}\mathcal{T}v_\alpha) \nonumber\\
&=&\frac{c_1+c_2}{16}\Big\{(l-1)g(\mathcal{T}v_\alpha,\mathcal{T}v_\alpha)
+3g(\mathcal{P}\mathcal{T}v_\alpha,\mathcal{P}\mathcal{T}v_\alpha)+2tr(f)g(f\mathcal{T}v_\alpha,\mathcal{T}v_\alpha)  \nonumber\\
&&-g(f\mathcal{T}v_\alpha,f\mathcal{T}v_\alpha)
+3g((f\mathcal{P}+t\mathcal{F})\mathcal{T}v_\alpha,(f\mathcal{P}+t\mathcal{F})\mathcal{T}v_\alpha)\Big\}  \\
&&+\frac{c_1-c_2}{16}\Big\{(l-2)g(f\mathcal{T}v_\alpha,\mathcal{T}v_\alpha)+tr(f)\cdot g(\mathcal{T}v_\alpha,\mathcal{T}v_\alpha) \nonumber\\
&&+6g(\mathcal{P}\mathcal{T}v_\alpha,(f\mathcal{P}+t\mathcal{F})\mathcal{T}v_\alpha) \Big\}  \nonumber\\
&&+\frac{1}{2}|[\mathcal{P},A_\alpha]|^2+tr(A_{\mathcal{N}v_\alpha}^2)-2tr(A_\alpha A_{\mathcal{N}v_\alpha}\mathcal{P})-tr(A_\alpha^2)  \nonumber\\
&=&\frac{c_1+c_2}{16}\Big\{-(l-1)\cdot tr(\mathcal{F}\mathcal{T})+3tr(\mathcal{P}^2\mathcal{T}\mathcal{F})-2tr(f)\cdot tr(f\mathcal{T}\mathcal{F})+tr(f^2\mathcal{T}\mathcal{F}) \nonumber\\
&&+3tr((f\mathcal{P}+t\mathcal{F})^2\mathcal{T}\mathcal{F}) \Big\}  \nonumber\\
&&+\frac{c_1-c_2}{16}\Big\{-(l-2)tr(f\mathcal{T}\mathcal{F})-tr(f)\cdot tr(\mathcal{F}\mathcal{T})+6tr((f\mathcal{P}+t\mathcal{F})\mathcal{T}\mathcal{F}\mathcal{P}) \Big\}  \nonumber\\
&&+\frac{1}{2}|[\mathcal{P},A_\alpha]|^2+tr(A_{\mathcal{N}v_\alpha}^2)-2tr(A_\alpha A_{\mathcal{N}v_\alpha}\mathcal{P})-tr(A_\alpha^2)  \nonumber
\end{eqnarray}

\textbf{Lemma 3.3} Let $M$ be an $l$-dimensional minimal submanifold in $\overline{M}^m\times \overline{M}^n$.
If $U$ is a parallel section in the normal bundle of $M$, then
\begin{eqnarray}
&&div(\nabla_{tU}tU) \\
&=&-div(tr(A_U f)tU)+\frac{c_1+c_2}{16}\Big\{(l-1)g(tU,tU)+3g(\mathcal{P}tU,\mathcal{P}tU) \nonumber\\
&&+2tr(f)g(ftU,tU)-g(ftU,ftU)+3g((f\mathcal{P}+t\mathcal{F})tU,(f\mathcal{P}+t\mathcal{F})tU)\Big\}   \nonumber\\
&&+\frac{c_1-c_2}{16}\Big\{(l-2)g(ftU,tU)+tr(f)\cdot g(tU,tU)
+6g(\mathcal{P}tU,(f\mathcal{P}+t\mathcal{F})tU) \Big\}  \nonumber\\
&&+tr(A_{sU}^2)-tr(A_U^2)+tr(A_U^2 f^2)+tr((A_U f)^2)-(tr(A_U f))^2   \nonumber\\
&&-2tr(A_U A_{sU}f)+\sum_\alpha g(A_U tv_\alpha,A_U tv_\alpha)-\sum_\alpha g(A_\alpha tU,A_\alpha tU)\nonumber
\end{eqnarray}
\textbf{Proof} As in the proof of Lemma 3.2, we have
\begin{eqnarray}
&&S(tU,tU) \nonumber\\
&=&\frac{c_1+c_2}{16}\Big\{(l-1)g(tU,tU)
+3g(\mathcal{P}tU,\mathcal{P}tU)+2tr(f)g(ftU,tU)  \nonumber\\
&&-g(ftU,ftU)+3g((f\mathcal{P}+t\mathcal{F})tU,(f\mathcal{P}+t\mathcal{F})tU)\Big\}  \nonumber\\
&&+\frac{c_1-c_2}{16}\Big\{(l-2)g(ftU,tU)+tr(f)\cdot g(tU,tU)
+6g(\mathcal{P}tU,(f\mathcal{P}+t\mathcal{F})tU) \Big\}  \nonumber\\
&&-\sum_\alpha g(A_\alpha tU,A_\alpha tU).
\end{eqnarray}

\begin{eqnarray}
\nabla_X(tU)=A_{sU}X-fA_U X
\end{eqnarray}
\begin{eqnarray}
div(tU)=-tr(A_U f)
\end{eqnarray}
\begin{eqnarray}
|\nabla tU|^2=tr(A_{sU}^2)+tr(A_U^2)-2tr(A_U A_{sU}f)-\sum_\alpha g(A_U tv_\alpha,A_U tv_\alpha )
\end{eqnarray}
\begin{eqnarray}
|L_{tU}g|^2=4tr(A_{sU}^2)+2tr(A_U^2 f^2)+2tr((A_U f)^2)-8tr(A_U A_{sU}f)
\end{eqnarray}
Put (3.31), (3.33)-(3.35) into
\begin{eqnarray}
div(\nabla_{tU}tU)&=&div(div(tU)tU)+S(tU,tU)  \nonumber\\
&&+\frac{1}{2}|L_{tU}g|^2-|\nabla tU|^2-(div(tU))^2 \nonumber
\end{eqnarray}
we get the result. $\Box$\vspace{2mm}

As before we choose an orthonormal basis $v_\alpha$ of $T(M)^\bot$ such that
$Dv_\alpha=0$ for all $\alpha$. Thus
\begin{eqnarray}
&&div(\nabla_{tv_\alpha}tv_\alpha) \\
&=&-div(tr(A_\alpha f)tv_\alpha)+\frac{c_1+c_2}{16}\Big\{(l-1)g(tv_\alpha,tv_\alpha)+3g(\mathcal{P}tv_\alpha,\mathcal{P}tv_\alpha) \nonumber\\
&&+2tr(f)g(ftv_\alpha,tv_\alpha)-g(ftv_\alpha,ftv_\alpha)
+3g((f\mathcal{P}+t\mathcal{F})tv_\alpha,(f\mathcal{P}+t\mathcal{F})tv_\alpha)\Big\}   \nonumber\\
&&+\frac{c_1-c_2}{16}\Big\{(l-2)g(ftv_\alpha,tv_\alpha)+tr(f)\cdot g(tv_\alpha,tv_\alpha)
+6g(\mathcal{P}tv_\alpha,(f\mathcal{P}+t\mathcal{F})tv_\alpha) \Big\}  \nonumber\\
&&+tr(A_{sv_\alpha}^2)-tr(A_{\alpha}^2)+tr(A_{\alpha}^2 f^2)+tr((A_{\alpha} f)^2)-(tr(A_{\alpha} f))^2
-2tr(A_{v_\alpha} A_{sv_\alpha}f) \nonumber\\
&=&-div(tr(A_\alpha f)tv_\alpha)+\frac{c_1+c_2}{16}\Big\{(l-1)\cdot tr(ht)-3tr(\mathcal{P}^2th)+2tr(f)\cdot tr(hft)  \nonumber\\
&&-tr(f^2th)-3tr((f\mathcal{P}+t\mathcal{F})^2th)\Big\}   \nonumber\\
&&+\frac{c_1-c_2}{16}\Big\{(l-2)\cdot tr(hft)+tr(f)\cdot tr(ht)-6tr((f\mathcal{P}+t\mathcal{F})th\mathcal{P})   \Big\}  \nonumber\\
&&+tr(A_{sv_\alpha}^2)-tr(A_{\alpha}^2)+tr(A_{\alpha}^2 f^2)+tr((A_{\alpha} f)^2)-(tr(A_{\alpha} f))^2
-2tr(A_{v_\alpha} A_{sv_\alpha}f) \nonumber
\end{eqnarray}

\textbf{Theorem 3.4}~~ If the submanifold $M$ is compact minimal submanifold in $\overline{M}^m\times \overline{M}^n$, then
\begin{eqnarray}
&&\int_M|\nabla A|^2d\mu+3(\frac{c_1+c_2}{16})^2\int_M W_1 d\mu+6(\frac{c_1-c_2}{16})^2\int_M W_2 d\mu  \nonumber \\
&=&-\int_M \sum_{i,j,\alpha}g((R(e_i,e_j)A_\alpha)e_i,A_\alpha e_j)d\mu\\
&&+3(\frac{c_1+c_2}{16})^2\int_M W_3 d\mu+6(\frac{c_1-c_2}{16})^2\int_M W_4 d\mu+\frac{c_1^2-c_2^2}{16^2}\int_M W_5 d\mu  \nonumber\\
&&+\frac{c_1+c_2}{16}\int_M W_6 d\mu+\frac{c_1-c_2}{16}\int_M W_7 d\mu  \nonumber
\end{eqnarray}
where
\begin{eqnarray}
W_1 &=&-tr(f^2\mathcal{T}\mathcal{F})+[tr(\mathcal{F}(\mathcal{P}t+\mathcal{T}s))]^2+2[tr(h(\mathcal{P}t+\mathcal{T}s))]^2  \\
&&+[tr((h\mathcal{P}+s\mathcal{F})(\mathcal{P}t+\mathcal{T}s))]^2  d\mu+[tr(\mathcal{F}\mathcal{T})]^2+[tr(h\mathcal{P}+s\mathcal{F})\mathcal{T}]^2 \nonumber\\
&&+2[tr(ht)]^2+[tr((h\mathcal{P}+s\mathcal{F})t)]^2+2tr((f\mathcal{P}+t\mathcal{F})^2)\cdot tr((\mathcal{F}t+\mathcal{N}s)^2) \nonumber
\end{eqnarray}
\begin{eqnarray}
W_2 &=&[tr(\mathcal{F}(\mathcal{P}t+\mathcal{T}s))]^2+tr(\mathcal{F}\mathcal{T})\cdot tr((h\mathcal{P}+s\mathcal{F})(\mathcal{P}t+\mathcal{T}s)) \nonumber\\
&&+tr((f\mathcal{P}+t\mathcal{F})^2)\cdot tr(\mathcal{N}^2)+tr(\mathcal{P}^2)\cdot tr((\mathcal{F}t+\mathcal{N}s)^2)
\end{eqnarray}
\begin{eqnarray}
W_3 &=&-(l-1)\cdot tr(\mathcal{F}\mathcal{T})+3tr(\mathcal{P}^2\mathcal{T}\mathcal{F})-2tr(f)\cdot tr(f\mathcal{T}\mathcal{F})  \nonumber\\
&&+3tr((f\mathcal{P}+t\mathcal{F})^2\mathcal{T}\mathcal{F})+tr[(\mathcal{F}(\mathcal{P}t+\mathcal{T}s))^2]  \nonumber\\
&&-2tr(\mathcal{N}(h\mathcal{P}+s\mathcal{F})\mathcal{P}(\mathcal{P}t+\mathcal{T}s))
+tr[(h(\mathcal{P}t+\mathcal{T}s))^2] \nonumber\\
&&+tr[((h\mathcal{P}+s\mathcal{F})(\mathcal{P}t+\mathcal{T}s))^2]  \\
&&-2tr[(\mathcal{F}t+\mathcal{N}s)(h\mathcal{P}+s\mathcal{F})(f\mathcal{P}+t\mathcal{F})(\mathcal{P}t+\mathcal{T}s)] \nonumber\\
&&+tr[(\mathcal{F}\mathcal{T})^2]
-2tr(\mathcal{N}\mathcal{F}\mathcal{P}\mathcal{T})+tr[(\mathcal{F}t)^2]  \nonumber\\
&&+tr[(h\mathcal{P}+s\mathcal{F})\mathcal{T}]^2
-2tr[(\mathcal{F}t+\mathcal{N}s)\mathcal{F}(f\mathcal{P}+t\mathcal{F})\mathcal{T}] \nonumber\\
&&+tr((\mathcal{F}t)^2)+2tr(\mathcal{N}h\mathcal{P}t)+tr((ht)^2)\nonumber\\
&&+tr(((h\mathcal{P}+s\mathcal{F})t)^2)+2tr((\mathcal{F}t+\mathcal{N}s)h(f\mathcal{P}+t\mathcal{F})t)  \nonumber\\
&&-2tr((f\mathcal{P}+t\mathcal{F})\mathcal{T}(\mathcal{F}t+\mathcal{N}s)\mathcal{F}) \nonumber\\
&&-2tr(\mathcal{P}(f\mathcal{P}+t\mathcal{F}))tr(\mathcal{N}(\mathcal{F}t+\mathcal{N}s))
+3tr(t(\mathcal{F}t+\mathcal{N}s)h(f\mathcal{P}+t\mathcal{F})) \nonumber\\
&&-2tr((\mathcal{P}t+\mathcal{T}s)(\mathcal{F}t+\mathcal{N}s)(h\mathcal{P}+s\mathcal{F})(f\mathcal{P}+t\mathcal{F})) \nonumber
\end{eqnarray}
\begin{eqnarray}
W_4 &=&-\Big[ -tr(\mathcal{T}\mathcal{F}(\mathcal{P}t+\mathcal{T}s)(h\mathcal{P}+s\mathcal{F}))
-tr(\mathcal{F}\mathcal{T}(h\mathcal{P}+s\mathcal{F})(\mathcal{P}t+\mathcal{T}s))
\nonumber \\
&&+2tr(\mathcal{T}\mathcal{N}(h\mathcal{P}+s\mathcal{F})(f\mathcal{P}+t\mathcal{F}))
+2tr((\mathcal{F}t+\mathcal{N}s)(h\mathcal{P}+s\mathcal{F})\mathcal{P}\mathcal{T})  \Big] \nonumber\\
&&-\Big[ tr(\mathcal{T}\mathcal{N}(h\mathcal{P}+s\mathcal{F})(f\mathcal{P}+t\mathcal{F}))
+tr(\mathcal{N}\mathcal{F}(f\mathcal{P}+t\mathcal{F})(\mathcal{P}t+\mathcal{T}s))  \nonumber \\
&&+tr(\mathcal{P}(f\mathcal{P}+t\mathcal{F}))\cdot tr(\mathcal{N}(\mathcal{F}t+\mathcal{N}s))  \Big] \nonumber\\
&&-\Big[ 2tr(\mathcal{P}\mathcal{T}(\mathcal{F}t+\mathcal{N}s)(h\mathcal{P}+s\mathcal{F}))  \\
&&+tr((f\mathcal{P}+t\mathcal{F})\mathcal{P})\cdot tr(\mathcal{N}(\mathcal{F}t+\mathcal{N}s))   \Big] \nonumber
\end{eqnarray}
\begin{eqnarray}
W_6 &=& -\frac{3}{2}\sum_{\alpha}|[f\mathcal{P}+t\mathcal{F}, A_\alpha]|^2-3\sum_{\alpha,\beta}g(A_\alpha tv_\beta, A_\beta tv_\alpha)  \nonumber\\
&&-\sum_{\alpha}\Big[2tr(A_\alpha^2)-2tr(A_{sv_\alpha}^2)+tr(A_\alpha A_{sv_\alpha}f)-tr(f)tr(A_\alpha A_{sv_\alpha}) \nonumber\\
&&+tr(f)tr(A_\alpha^2f)-2tr(A_\alpha^2f^2)-2(tr(A_\alpha f))^2+tr((A_\alpha f)^2) \nonumber\\
&&+3tr(A_\alpha A_{(h\mathcal{P}+s\mathcal{F})(\mathcal{P}t+\mathcal{T}s)v_\alpha})\Big]
\end{eqnarray}
\begin{eqnarray}
W_5 &=& 3\Big[-(l-2)tr(f\mathcal{T}\mathcal{F})-tr(f)\cdot tr(\mathcal{F}\mathcal{T})
+6tr((f\mathcal{P}+t\mathcal{F})\mathcal{T}\mathcal{F}\mathcal{P})\Big] \nonumber\\
&&-6\Big[tr(\mathcal{F}(\mathcal{P}t+\mathcal{T}s))\cdot tr((h\mathcal{P}+s\mathcal{F})(\mathcal{P}t+\mathcal{T}s))  \nonumber\\
&&-tr(\mathcal{F}(\mathcal{P}t+\mathcal{T}s)(h\mathcal{P}+s\mathcal{F})(\mathcal{P}t+\mathcal{T}s)) \nonumber\\
&&+tr(\mathcal{N}(h\mathcal{P}+s\mathcal{F})(f\mathcal{P}+t\mathcal{F})(\mathcal{P}t+\mathcal{T}s))   \\
&&+tr((\mathcal{F}t+\mathcal{N}s)(h\mathcal{P}+s\mathcal{F})\mathcal{P}(\mathcal{P}t+\mathcal{T}s))\Big] \nonumber\\
&&-6\Big[tr(\mathcal{F}\mathcal{T})\cdot tr((h\mathcal{P}+s\mathcal{F})\mathcal{T})
-tr(\mathcal{F}\mathcal{T}(h\mathcal{P}+s\mathcal{F})\mathcal{T})  \nonumber\\
&&+tr(\mathcal{N}\mathcal{F}(f\mathcal{P}+t\mathcal{F})\mathcal{T})+tr(\mathcal{F}\mathcal{P}\mathcal{T}(\mathcal{F}t+\mathcal{N}s))\Big] \nonumber\\
&&-6\Big[ -tr((h\mathcal{P}+s\mathcal{F})t\mathcal{F}t)
-tr((f\mathcal{P}+t\mathcal{F})t\mathcal{N}h)-tr(h\mathcal{P}t(\mathcal{F}t+\mathcal{N}s)) \Big] \nonumber\\
&&-6\Big[tr(\mathcal{T}(\mathcal{F}t+\mathcal{N}s)(h\mathcal{P}+s\mathcal{F})(f\mathcal{P}+t\mathcal{F})) \nonumber\\
&&+tr((\mathcal{P}t+\mathcal{T}s)(\mathcal{F}t+\mathcal{N}s)\mathcal{F}(f\mathcal{P}+t\mathcal{F})) \nonumber\\
&&+tr((f\mathcal{P}+t\mathcal{F})^2)\cdot tr(\mathcal{N}(\mathcal{F}t+\mathcal{N}s)) \nonumber\\
&&+tr((\mathcal{F}t+\mathcal{N}s)^2)\cdot tr(\mathcal{P}(f\mathcal{P}+t\mathcal{F})) \Big] \nonumber\\
&&-6\Big[ tr(\mathcal{F}\mathcal{T})\cdot tr(\mathcal{F}(\mathcal{P}t+\mathcal{T}s))
-tr(\mathcal{F}\mathcal{T}\mathcal{F}(\mathcal{P}t+\mathcal{T}s)) \nonumber \\
&&+2tr(\mathcal{N}\mathcal{F}\mathcal{P}(\mathcal{P}t+\mathcal{T}s))+2tr(h\mathcal{T})\cdot tr(h(\mathcal{P}t+\mathcal{T}s)) \nonumber \\
&&-tr(h\mathcal{T}h(\mathcal{P}t+\mathcal{T}s))
+tr((h\mathcal{P}+s\mathcal{F})\mathcal{T})\cdot tr((h\mathcal{P}+s\mathcal{F})(\mathcal{P}t+\mathcal{T}s)) \nonumber \\
&&-tr((h\mathcal{P}+s\mathcal{F})\mathcal{T}(h\mathcal{P}+s\mathcal{F})(\mathcal{P}t+\mathcal{T}s))  \nonumber\\
&&+2tr((\mathcal{F}t+\mathcal{N}s)(h\mathcal{P}+s\mathcal{F})(f\mathcal{P}+t\mathcal{F})\mathcal{T})  \Big] \nonumber\\
&&+3\Big[ -2tr(\mathcal{T}\mathcal{N}\mathcal{F}(f\mathcal{P}+t\mathcal{F}))
-2tr(\mathcal{P}(f\mathcal{P}+t\mathcal{F}))\cdot tr(\mathcal{N}^2) \nonumber \\
&&+3tr(t\mathcal{N}h(f\mathcal{P}+t\mathcal{F}))-2tr((\mathcal{P}t+\mathcal{T}s)\mathcal{N}(h\mathcal{P}+s\mathcal{F})(f\mathcal{P}+t\mathcal{F})) \nonumber \\
&&-2tr((f\mathcal{P}+t\mathcal{F})^2)\cdot tr(\mathcal{N}(\mathcal{F}t+\mathcal{N}s)) \Big] \nonumber\\
&&+3\Big[ -2tr(\mathcal{F}\mathcal{P}\mathcal{T}(\mathcal{F}t+\mathcal{N}s))
-2tr(\mathcal{P}^2)tr(\mathcal{N}(\mathcal{F}t+\mathcal{N}s))  \nonumber\\
&&+3tr(h\mathcal{P}t(\mathcal{F}t+\mathcal{N}s))
-2tr((\mathcal{P}t+\mathcal{T}s)(\mathcal{F}t+\mathcal{N}s)(h\mathcal{P}+s\mathcal{F})\mathcal{P})  \nonumber\\
&&-2tr((f\mathcal{P}+t\mathcal{F})\mathcal{P})\cdot tr((\mathcal{F}t+\mathcal{N}s)^2) \Big] \nonumber
\end{eqnarray}
\begin{eqnarray}
W_7 &=& -l\cdot \Big[tr(A_\alpha^2f)-tr(A_\alpha A_{sv_\alpha})\Big]
-6\Big[ tr(A_\alpha A_{(h\mathcal{P}+s\mathcal{F})\mathcal{T}v_\alpha}) \nonumber\\
&&-tr(A_\alpha^2(f\mathcal{P}+t\mathcal{F})\mathcal{P})
+tr(A_\alpha\mathcal{P}A_\alpha(f\mathcal{P}+t\mathcal{F}))  \Big]
\end{eqnarray}

Note that all the terms in $W_1$ and $W_2$ are positive,  the left side of $(3.39)$ is positive, therefore
we get the following

\textbf{Corollary 3.5}  Let $M$ be an $l$-dimensional compact minimal submanifold
in $\overline{M}^m\times \overline{M}^n$
with flat normal connection. If the sectional curvature K of M has a
lower bound $C$ and if the second fundamental form $A$ satisfies
\begin{eqnarray}
&&-l\cdot C\sum_\alpha tr(A_\alpha^2)+3(\frac{c_1+c_2}{16})^2 W_3+6(\frac{c_1-c_2}{16})^2 W_4\nonumber\\
&&+\frac{c_1^2-c_2^2}{16^2} W_5+\frac{c_1+c_2}{16}W_6+\frac{c_1-c_2}{16}W_7<0
\end{eqnarray}
then $M$ is totally geodesic.

\textbf{Proof}. We can choosing an orthonormal basis $\{e_i\}$ of $T_x M$ such that $A_\alpha e_i=h_i^\alpha e_i,$ $i=1,2,\cdots$, then
\begin{eqnarray}
&&\sum_{i,j}((R(e_i,e_j)A)_\alpha e_i,A_\alpha e_j)\\
&=&\sum_{i,j}g(R(e_i,e_j)A_\alpha e_i,A_\alpha e_j)-\sum_{i,j}g(A_\alpha R(e_i,e_j)e_i,A_\alpha e_j) \nonumber\\
&=&\frac{1}{2}\sum_{i,j}(h_i^\alpha-h_j^\alpha)^2K_{ij}  \nonumber\\
&\geq & \frac{C}{2}\sum_{i,j}(h_i^\alpha-h_j^\alpha)^2=l\cdot C\cdot tr(A_\alpha^2)   \nonumber
\end{eqnarray}
where $K_{ij}$ denotes the sectional curvature of $M$ with respect to the section spanned by $e_i$ and $e_j$.

Denote the left side of (3.39) by $\mathcal{L}$, then there exist a small $\varepsilon >0$, such that
\begin{eqnarray}
&&\mathcal{L}+\varepsilon\int_M\sum_\alpha tr(A_\alpha^2)d\mu \nonumber\\
&&\leq \int_M \Big[(\varepsilon-l\cdot C)\sum_\alpha tr(A_\alpha^2)+3(\frac{c_1+c_2}{16})^2 W_3+6(\frac{c_1-c_2}{16})^2 W_4\nonumber\\
&&+\frac{c_1^2-c_2^2}{16^2} W_5+\frac{c_1+c_2}{16}W_6+\frac{c_1-c_2}{16}W_7\Big] d\mu \nonumber\\
&&\leq 0
\end{eqnarray}
Therefore $\sum_\alpha tr(A_\alpha^2)=0$. $\Box$\vspace{2mm}

\textbf{Corollary 3.6}~~Let $M$ be an $l$-dimensional compact minimal submanifold
in $\overline{M}^m\times \overline{M}^n$
with flat normal connection, $c_1+c_2>0$. If the sectional curvature $K$ of $M$ has a
lower bound $C$ and if the second fundamental form $A$ satisfies
\begin{eqnarray}
&&-l\cdot C\sum_\alpha tr(A_\alpha^2)+3(\frac{c_1+c_2}{16})^2 W_3+6(\frac{c_1-c_2}{16})^2 W_4 \\
&&+\frac{c_1^2-c_2^2}{16^2} W_5+\frac{c_1+c_2}{16}W_6'+\frac{c_1-c_2}{16}W_7\leq 0,  \nonumber
\end{eqnarray}
where $W_6'=W_6+2\sum_\alpha tr(A_\alpha^2)+\frac{3}{2}\sum_{\alpha}|[f\mathcal{P}+t\mathcal{F}, A_\alpha]|^2$,
then $M$ is totally geodesic.

\textbf{Proof}. ~~We only need to move the terms $-2\sum_\alpha tr(A_\alpha^2)$ and \\
 $-\frac{3}{2}\sum_{\alpha}|[f\mathcal{P}+t\mathcal{F}, A_\alpha]|^2$ on the left side of $(3.39)$.
The rest are the same as the proof of Corollary 3.5.
$\Box$\vspace{2mm}

\section{The $F$-anti invariant submanifold }
From now on we assume that $M$ is $F$-anti invariant, that is $FT_x M\subset T_x(M)^\bot$, or equivalently $f\equiv 0$.
Then from Lemma 3.2 or (3.31) we have
\begin{eqnarray}
&&\sum_\alpha\Big[\frac{1}{2}|[\mathcal{P},A_\alpha]|^2+tr(A_{\mathcal{N}v_\alpha}^2)-2tr(A_\alpha A_{\mathcal{N}v_\alpha}\mathcal{P})-tr(A_\alpha^2)\Big] \nonumber\\
&=&-\frac{c_1+c_2}{16}\sum_\alpha\Big[(l-1)g(\mathcal{T}v_\alpha,\mathcal{T}v_\alpha)
+3g(\mathcal{P}\mathcal{T}v_\alpha,\mathcal{P}\mathcal{T}v_\alpha)
+3g(t\mathcal{F}\mathcal{T}v_\alpha,t\mathcal{F}\mathcal{T}v_\alpha)\Big] \nonumber\\
&&-6\frac{c_1-c_2}{16}\sum_\alpha g(\mathcal{P}\mathcal{T}v_\alpha, t\mathcal{F}\mathcal{T}v_\alpha)+\sum_\alpha div(\nabla_{\mathcal{T}v_\alpha}\mathcal{T}v_\alpha)\\
&=&-\frac{c_1+c_2}{16}\Big[-(l-1)tr(\mathcal{F}\mathcal{T})+3tr(\mathcal{P}^2\mathcal{T}\mathcal{F})+3tr(ht(\mathcal{F}\mathcal{T})^2)\Big] \nonumber\\
&&-6\frac{c_1-c_2}{16}tr(t\mathcal{F}\mathcal{T}\mathcal{F}\mathcal{P})+\sum_\alpha div(\nabla_{\mathcal{T}v_\alpha}\mathcal{T}v_\alpha)\nonumber
\end{eqnarray}
and from Lemma 3.3 or (3.38) we have
\begin{eqnarray}
&&\sum_\alpha\Big[tr(A_{\alpha}^2)-tr(A_{sv_\alpha}^2)\Big] \nonumber\\
&=&\frac{c_1+c_2}{16}\sum_\alpha\Big[(l-1)g(tv_\alpha,tv_\alpha)
+3g(\mathcal{P}tv_\alpha,\mathcal{P}tv_\alpha)+3g(t\mathcal{F}tv_\alpha,t\mathcal{F}tv_\alpha)\Big] \nonumber\\
&&+6\frac{c_1-c_2}{16}\sum_\alpha g(\mathcal{P}tv_\alpha, t\mathcal{F}tv_\alpha)-\sum_\alpha div(\nabla_{tv_\alpha}tv_\alpha) \\
&=&\frac{c_1+c_2}{16}\Big[(l-1)tr(ht)-3tr(\mathcal{P}^2th)-3tr(h\mathcal{T}ht\mathcal{F}t)\Big] \nonumber\\
&&-6\frac{c_1-c_2}{16}tr(t\mathcal{F}th\mathcal{P})-\sum_\alpha div(\nabla_{tv_\alpha}tv_\alpha)\nonumber
\end{eqnarray}

\textbf{Theorem 4.1} Let $M$ be an $F$-anti invariant submanifold in an arbitrary kaehler product manifold
$\overline{M}^m\times\overline{M}^n$, if $s\equiv 0$, then $M$ is totally geodesic.

\textbf{Proof}.  From \cite{YK1} Theorem 3.3 on Page 426, we know that when
$M$ is $F$-anti invariant,
$A_{hX}Y=0$ or equivalently $tB(X,Y)=0$
for all tangent vector fields $X, Y$ of $M$.
Then by (3.14)
$$\sum_\alpha\Big[tr(A_{\alpha}^2)-tr(A_{sv_\alpha}^2)\Big]=0.$$

Therefore   $\sum_\alpha tr(A_{\alpha}^2)=\sum_\alpha tr(A_{sv_\alpha}^2)=0.$  $\Box$\vspace{2mm}

Similarly
$$\sum_{\alpha,\beta}g(A_\alpha tv_\beta, A_\beta tv_\alpha)=\sum_{i,\alpha}g(e_i,A_{hA_\alpha e_i}tv_\alpha)=0. $$

\textbf{Theorem 4.2} Let $M$ be an $F$-anti invariant Lagrangian submanifold in an arbitrary kaehler product manifold
$\overline{M}^m\times\overline{M}^n$, then $M$ is totally geodesic.

\textbf{Proof}. The dimension of Lagrangian submanifold is half of that of $\overline{M}^m\times\overline{M}^n$,
by Theorem 3.3 on Page 426 of \cite{YK1}, we get the result.   $\Box$\vspace{2mm}

\textbf{Remark 4.3} Certainly Corollary 3.5, Corollary 3.6 are still true in this case and
$W_1,\cdots , W_7 $ have simpler expressions, specially
\begin{eqnarray}
W_6 &=& -\frac{3}{2}\sum_{\alpha}|[t\mathcal{F}, A_\alpha]|^2
-3\sum_{\alpha}tr(A_\alpha A_{(h\mathcal{P}+s\mathcal{F})(\mathcal{P}t+\mathcal{T}s)v_\alpha}) \nonumber
\end{eqnarray}

\end{document}